\documentclass[12pt]{amsart}
\usepackage{amssymb,latexsym}
\usepackage{amscd}
\usepackage{verbatim}
\usepackage{amsmath,amsfonts,amssymb,epsf}
\usepackage{graphicx}
\newcommand{\R}{\mathbb{R}}
\newcommand{\eps}{\varepsilon}
\newtheorem{theorem}{Theorem}[section]
\newtheorem{lemma}[theorem]{Lemma}

\newtheorem{remark}[theorem]{Remark}
\newtheorem{assumption}[theorem]{Assumptions}
\newtheorem{corollary}[theorem]{Corollary}

\newtheorem{example}[theorem]{Example}
\newtheorem{definition}[theorem]{Definition}

\newcommand{\dxx}{\,\mathrm{d}\mathbf{x}}

\newcommand{\dt}{\,\mathrm{d}t}

\newcommand{\dss}{\,\mathrm{d}\mathbf{s}}

\newcommand{\dtet}{\,\mathrm{d}\theta}
\newcommand{\be}{\begin{equation}}
\newcommand{\ee}{\end{equation}}
\newcommand{\mysection}[1]{\section{#1}\setcounter{equation}{0}}
\newcommand{\bea}{\begin{eqnarray}}
\newcommand{\eea}{\end{eqnarray}}
\newcommand{\bean}{\begin{eqnarray*}}
\newcommand{\eean}{\end{eqnarray*}}

\def\squarebox#1{\hbox to #1{\hfill\vbox to #1{\vfill}}}

\begin{document}
\renewcommand{\theequation}{\thesection.\arabic{equation}}
\title[Schr\"{o}dinger operators on
infinite trees]%
{Spectral properties of Schr\"{o}dinger operators defined on
$N$-dimensional  infinite trees}
\author{Yehuda Pinchover}
\address{Department of Mathematics\\ Technion - Israel Institute of
Technology\\ Haifa 32000, Israel}
\email{pincho@techunix.technion.ac.il}
\author{Gershon Wolansky}
\address{Department of Mathematics\\ Technion - Israel Institute of
Technology\\ Haifa 32000, Israel}
 \email{gershonw@math.technion.ac.il}
\author{Daphne Zelig}
\address{Department of Mathematics\\ Technion - Israel Institute of
Technology\\ Haifa 32000, Israel}
 \email{zeligd@techunix.technion.ac.il}
 %
\date{June 6, 2006}
\begin{abstract}
We study the discreteness of the spectrum of Schr\"{o}\-dinger
operators which are defined on $N$-dimensional rooted trees of a
finite or infinite volume, and are subject to a certain mixed
boundary condition.  We present a method to estimate their
eigenvalues using operators on a one-dimensional tree. These
operators are called {\em width-weighted operators}, since their
coefficients depend on the section width or area of the
$N$-dimensional tree.  We show that the spectrum of the
width-weighted operator tends to the spectrum of a one-dimensional
limit operator as the sections width tends to zero. Moreover, the
projections to the one-dimensional tree of eigenfunctions of the
$N$-dimensional Laplace operator converge to the corresponding
eigenfunctions of the one-dimensional limit operator.
\end{abstract}
\maketitle
\begin{figure}[t] \centering
\includegraphics[scale=0.4]{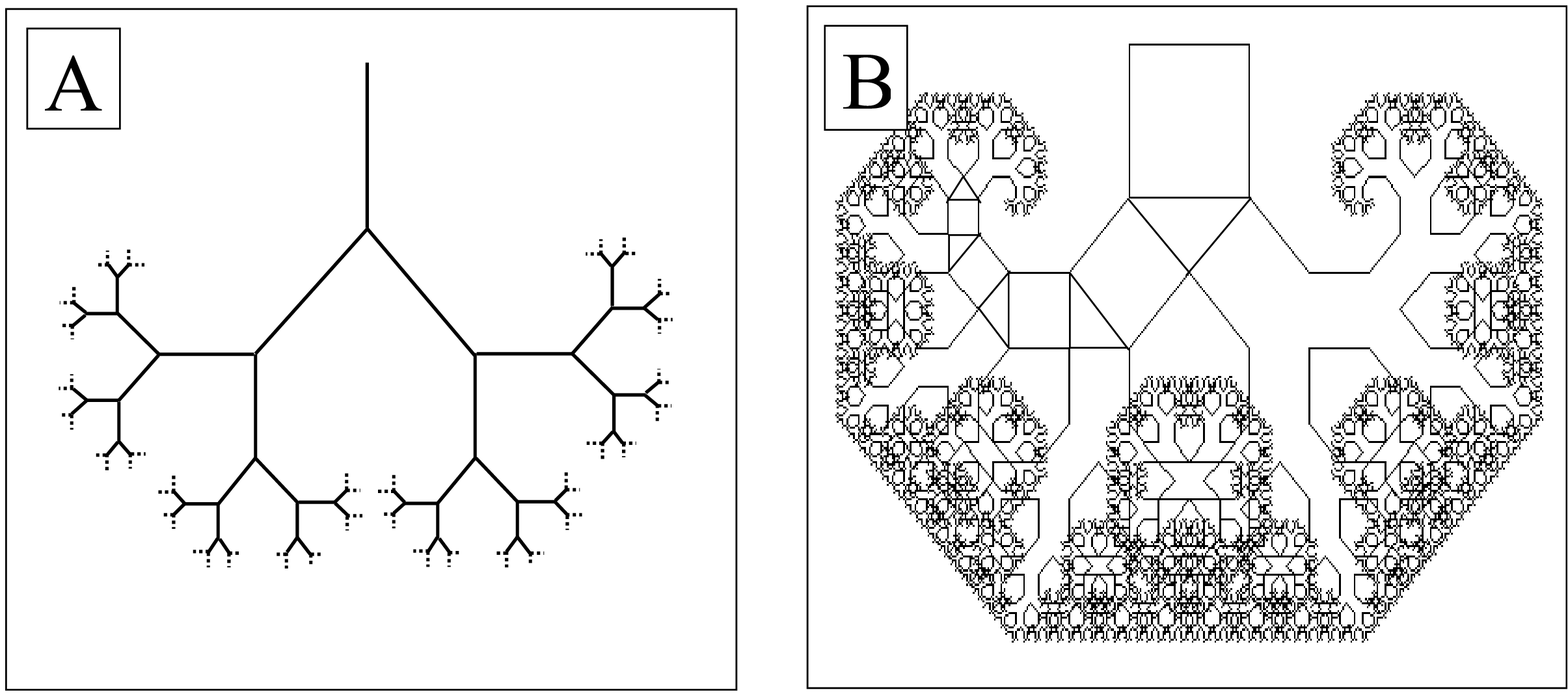}
\caption{An example of one and two dimensional trees.}{{\bf A.}
One-dimensional tree. {\bf B.} Two-dimensional tree presented in
$\mathbb{R}^2$. Some of its triangle connectors and rectangle
edges are emphasized.} \label{T2AndT1}
\end{figure}

\tableofcontents
\mysection{Introduction} Let $T_1$ be a one-dimensional infinite
tree. We assume throughout this paper that $T_1$ is regular (see
Definition~\ref{defregular} and Remark~\ref{remgenarlity}). For
$N\geq 2$, we also consider an $\eps$-inflated tree $T_N^\eps$
around $T_1$ which is an $N$-dimensional offset (or fattening) of
$T_1$. See Figure~\ref{T2AndT1} for an illustration
 of a $2$-dimensional tree $T_2=T^1_2$.

We prove $\eps$-dependent estimates for the spectrum of the
eigenvalue problem
\begin{equation}\label{SchrDef}
 L_\eps u:=-\Delta u+W_{T_N^\eps}u=\lambda^\eps u \qquad \mbox{  on  } H^1_0(T_N^\eps),
\end{equation}
subject to the Neumann boundary condition  on $\partial T^\eps_N$
except on the root of the tree, where we impose the Dirichlet
boundary condition. We assume that  $W_{T_N^\eps}$ is a bounded
and continuous potential on $T^\eps_N$.
Specifically, we show that if $T_N^\eps$ has a finite radius, then
under some further assumptions the spectrum of $L_\eps$ is
discrete and  the eigenvalues of the Schr\"odinger operators
$L_\eps$ satisfy $\lambda_i^\eps\rightarrow \mu_i$ as
$\eps\rightarrow 0$, where $\mu_i$ are the eigenvalues of the
following {\em weighted} Schr\"odinger operator on $T_1$
\begin{equation}\label{AnDef}
 \overline{L} u:=-\frac{1}{\rho}(\rho
u')'+W_{T_1}u = \mu_i u .
\end{equation}
Here $\rho>0$ is a weight function on $T_1$ defined in terms of
the inflation $T^\eps_N$, and $W_{T_1}$ is the cross section
average of $W_{T_N^\eps}$. \par The spectral behavior of the
Neumann Laplacian and Schr\"{o}dinger operators on thin domains
has been extensively investigated. Indeed, in \cite{Rub1},
Rubinstein and Schatzman study the relation between the spectral
properties of the Laplace operator defined on a metric {\em graph}
$G$ and on a strip shaped domain $G^\eps$ of width $\eps$ around
$G$.
The results of \cite{Rub1} on the spectrum of the Laplacian cannot
be applied to our trees because of the following essential
differences between the problems:
\begin{enumerate}
\item{Rubinstein and Schatzman treat the case in which the graph
$G$  has a {\em finite}
        number of vertices, while our tree $T_1$ has an infinite number of vertices.}
\item{They consider graph-surrounding domains having a constant
(uniform)
        width. In the case of an infinite trees, the
        discreteness of the spectrum imposes that
the width of higher branches of the tree must be scaled. }
       \item{In particular, the inflated finite graph is of
 finite volume, while our inflated infinite tree may have infinite volume.  }
\end{enumerate}
In \cite{Kuch2}, Kuchment and Zeng extend the results in
\cite{Rub1}. For example, the conditions on the smoothness of the
boundary of the domain near the vertices were relaxed and the
constant width of the surrounding domain is not assumed.

Since $T_1$ in our case is an infinite tree, the results of
\cite{Kuch2, Rub1} do not apply in our setting.  Nevertheless, we
were able to modify the approach in \cite{Rub1} to obtain similar
results in the infinite case. In particular, we could not compare
directly the eigenvalue $\lambda^\eps_i$ to $\mu_i$. Instead, we
find it more convenient to compare the spectra of
$-\Delta+W_{T_N^\eps}$ on $T^\eps_N$ to the Schr\"odinger operator
on $T_1$ subjected to a pair of $\eps-$dependent weight functions
$\rho_{1,\eps}$, $\rho_{2,\eps}$, satisfying $\rho_{1,\eps},
\rho_{2,\eps} \rightarrow \rho$ as $\eps\rightarrow 0$.  So, we
replace (\ref{AnDef}) by
$$ \overline{L}_\eps u := -\frac{1}{\rho_{2,\eps}}(\rho_{1,\eps}
u')'+W_{T_1}u  = \mu_i^\eps u \, , $$ and prove that the
$\lambda_j^\eps$ is approximated, on the one hand, by $\mu_j^\eps$
while the later is approximated by $\mu_j$ for $\eps$ small.

Spectral properties of Schr\"{o}dinger operators defined on
infinite one-dimensional metric trees and graphs has also been
intensively studied. In \cite{Car}, Carlson shows that if $G$ is a
connected metric graph which has a finite total edges length (a
finite volume), then the Laplacian defined on $G$ has a compact
resolvent and therefore a discrete spectrum. Solomyak and Naimark
have developed general tools for studying spectral properties of
Schr\"{o}dinger operators on metric graphs and trees (see, for
example, \cite{Naim1, Naim2, Sol1, Sol2}).
 In \cite{Sol1}, Solomyak has proved that if $T_1$ is a
regular tree whose radius is finite, and if $W_{T_1}(x)$ is a
radial measurable real valued function which is bounded below,
then the spectrum of $\overline{L}$ is discrete.

Solomyak's result is stated for trees of uniform weight function
$\rho$ and its proof relies on the monotonicity of $g$, where
$g(t)$ is the number of branches which contain points of distance
$t$ from the root. In fact, to adjust Solomyak's proof for our
case, one needs to assume only that $g\rho$ is a monotone
nondecreasing. If $\rho$ is constant then it is a natural
assumption, but if $\rho(t)$ is decreasing (as in our case), this
monotonicity may be violated. So, we extend this result under a
milder condition on $g\rho$.

We prove the discreteness of the spectrum of Schr\"{o}dinger
operators on regular $N$-dimensional trees with infinite volume,
as long as the tree radius is finite. Our proof relies on a lemma
of Lewis \cite[Lemma 1]{Lew}. The proof of the discreteness in the
$N$-dimensional case can be applied also to show that the
$L^2$-norm of functions which are bounded in $H^1_0(T_N^\eps)$
does not accumulate at the tree connectors or ends.

A natural question emerging from the correspondence between the
eigenvalues of $N$-dimensional Laplace operator, and
one-dimensional width-weighted operators, is whether the
corresponding {\it eigenfunctions} present the same convergence
behavior. In \cite{Kosugi1,Kosugi2}, Kosugi has proved that
solutions of (semilinear) elliptic equations on finite
$N$-dimensional trees indeed converge as the width tends to zero
to solutions  of width-weighted equations. We present a different
method and prove that certain projections of eigenfunctions of the
Laplace operator on $T_N^\eps$ converge to the corresponding
eigenfunctions on $T_1$. In contrast to \cite{Kosugi1,Kosugi2}, we
treat infinite trees rather than trees with a finite number of
vertices. In addition, our assumptions on the smoothness of the
connectors are much weaker than those in \cite{Kosugi1,Kosugi2},
and in fact, we require only that the connectors have a Lipschitz
boundary.

\begin{remark}\label{remgenarlity}{\em
 Our method applies to  more
general setting. But to facilitate the presentation,  we restrict
our study in the present paper to the case where $T_1$ is a
regular metric tree, and the inflated $N$-dimensional tree is a
self-similar radial tree with `cylindrical' edges.
 } \end{remark}

We wish to mention two more articles which study the spectrum of
thin domains. In an earlier article \cite{Kuch1}, Kuchment and
Zeng study the dependence of the spectrum of the Neumann Laplacian
on the behavior of the surrounding thin domain near the vertices.
They found differential operators on the graph which correspond to
the case in which the neighborhoods of the vertices are much
larger or much smaller than the tubes connecting them. In
\cite{Eva}, Evans and Saito proved results about the connection
between the essential spectrum of the Neumann Laplacian on thin
domains surrounding trees and the essential spectrum of their
skeletons. They apply their results on horns, spirals, ``rooms and
passages" domains and domains with fractal boundaries. In our case
 the essential spectrum is empty,  as was mentioned above.

The motivation for our problem is that fractal structures, and in
particular, fractal tree-like structures, have a vast applications
range. For example, fractal geometry is used in order to form
antennas, which present a multi-band behavior (see
\cite{Anag,Puent1}). In \cite{Puent2}, Puente et al. state that
fractal tree shaped antennas have a denser band distribution than
previously reported Sierpinski fractal antennas. Estimating the
eigenvalues of the Laplace operator defined on such domains may
help in specifying the natural transmission frequencies for the
antennas.

 Another applications field for fractal geometry is
medical modelling. Nelson et al. mention in \cite{Nel1} that
fractal models can be applied to human lungs, vascular tree,
neural networks, urinary ducts, brain folds and cardiac conduction
fibers. Fractal models of human lungs can be found also in
\cite{Man,Nel2,Weib}.

\vskip 4mm

The outline of this article is as follows. In
Section~\ref{SecNotations}, we present the basic notations we use,
describe the class of trees we are interested in, and define the
operators on the trees. Section~\ref{appen1} is devoted to the
study of the behavior of $H^1$-functions near the vertices. In
Section~\ref{SectionBackground}, we prove the discreteness of the
spectrum of Schr\"{o}dinger operators on $T_1$ and $T_N$. The
convergence (as $\eps\rightarrow 0$)  of the spectrum of
$\{\overline{L}_{\eps}\}$, the operator sequence defined on $T_1$,
to the spectrum of the limit operator $\overline{L}$ is proved in
Section~\ref{AnLimitSec}.

In sections~\ref{QSection} and \ref{PSection} we define
transformations between $H_0^1(T_N)$ and $H^1_{0,\rho_2}(T_1)$ and
prove comparison theorems for the Rayleigh quotients of the one
and $N$-dimensional operators. In Section~\ref{SectionKobyLemma},
we use these comparison theorems to characterize the behavior of
the spectrum on  $T_N$. Finally, the convergence of projections of
$N$-dimensional eigenfunctions of Laplace operator to
eigenfunctions of the one-dimensional width-weighted operators is
proved in Section~\ref{SecEigHarmonicConv}.

\mysection{Preliminaries} \label{SecNotations}
\subsection{General notations}
\begin{enumerate}
\item{Throughout the article, $c,c_1,c_2,\ldots$, and $C$ denote
constants, whose exact values are irrelevant, and may change from
line to line.}
\item{Let $\{a_j\}$ and $\{b_j\}$ be positive sequences. We denote
$a_j\asymp b_j$ if there exists a constant $c>0$ such that
$c^{-1}\leq {a_j}/{b_j}\leq c$ for all $j\in \mathbb{N}$. We use a
similar notation for positive functions, i.e., we denote $f\asymp
g$ if there exists a constant $c>0$ such that $c^{-1}\leq
{f(x)}/{g(x)}\leq c$ for all $x$ in the domain of the functions
$f$ and $g$.}
\item{For a domain $\Omega\subset \mathbb{R}^N$, we denote by
$|\Omega|$ its volume in $\mathbb{R}^N$.}
\end{enumerate}

\subsection{The tree $T_1$}\label{3.2}
\begin{enumerate}
 \item $T_1$ is a one-dimensional connected rooted metric tree. It
contains an infinite number of vertices $v$, connected by  edges
$e$.
 \item The root $O$ of $T_1$ is a distinguished (and unique)
vertex. Its generation number is defined to be zero.
 \item A vertex of
$T_1$ is of generation  $j$ if it is connected to the root by a
succession of $j$ edges. The generation of a given vertex $v$ is
denoted by $\mathrm{gen}(v)$.
 \item  Likewise, $e$ is an edge of
generation $j$ if it connect a pair of vertices of generations $j$
and $j+1$, respectively. The generation number of a given edge $e$
is denoted by $\mathrm{gen}(e)$.
 \item The Euclidian length of an edge $e$
is denoted by $|e|$.
 \item The degree of a vertex $v$ is $k(v)$.
It is the number of edges connecting $v$ to the vertices of
generation $\mathrm{gen}(v)+1$.
 \item The set of all edges meeting
at a vertex $v$ is $N(v)$. There are exactly $k(v)+1$ edges in
$N(v)$.
 \item The distance $\mathrm{dist}(x,y)$ between $x,y\in T_1$ is
the Euclidian length of  the path on $T_1$  connecting $x$ to $y$.
We denote $|x|:={\rm dist}\{O,x\}$.
 \item $g(t)$ is the counting function of $T_1$, namely, $g(t)$ is
the number of edges which contain a point $x\in T_1$ with $|x|=t$.
 \item $R(T_1)\equiv\sup_{x\in T_1}|x|$ is
the {\em radius} of $T_1$. $L(T_1)\equiv \sum_{e\in T_1} |e|$ is
the {\em length} of $T_1$.
\end{enumerate}
\begin{definition}\label{defregular} {\em
 $T_1$ is called  {\em radial} if  the  length $|e|$ of each edge $e$ and
 the degree $k(v)$ of each vertex $v$ depend only on $\mathrm{gen}(e)$ and $\mathrm{gen}(v)$,
 respectively. A radial tree is called {\em regular} if $k(v)=k$ is a
 constant, independent of the generation.
} \end{definition}

\subsection{The $\eps$-inflated $N$-dimensional
tree}\label{epsinf} 
\begin{figure}[t]
\centering
\includegraphics[scale=0.65]{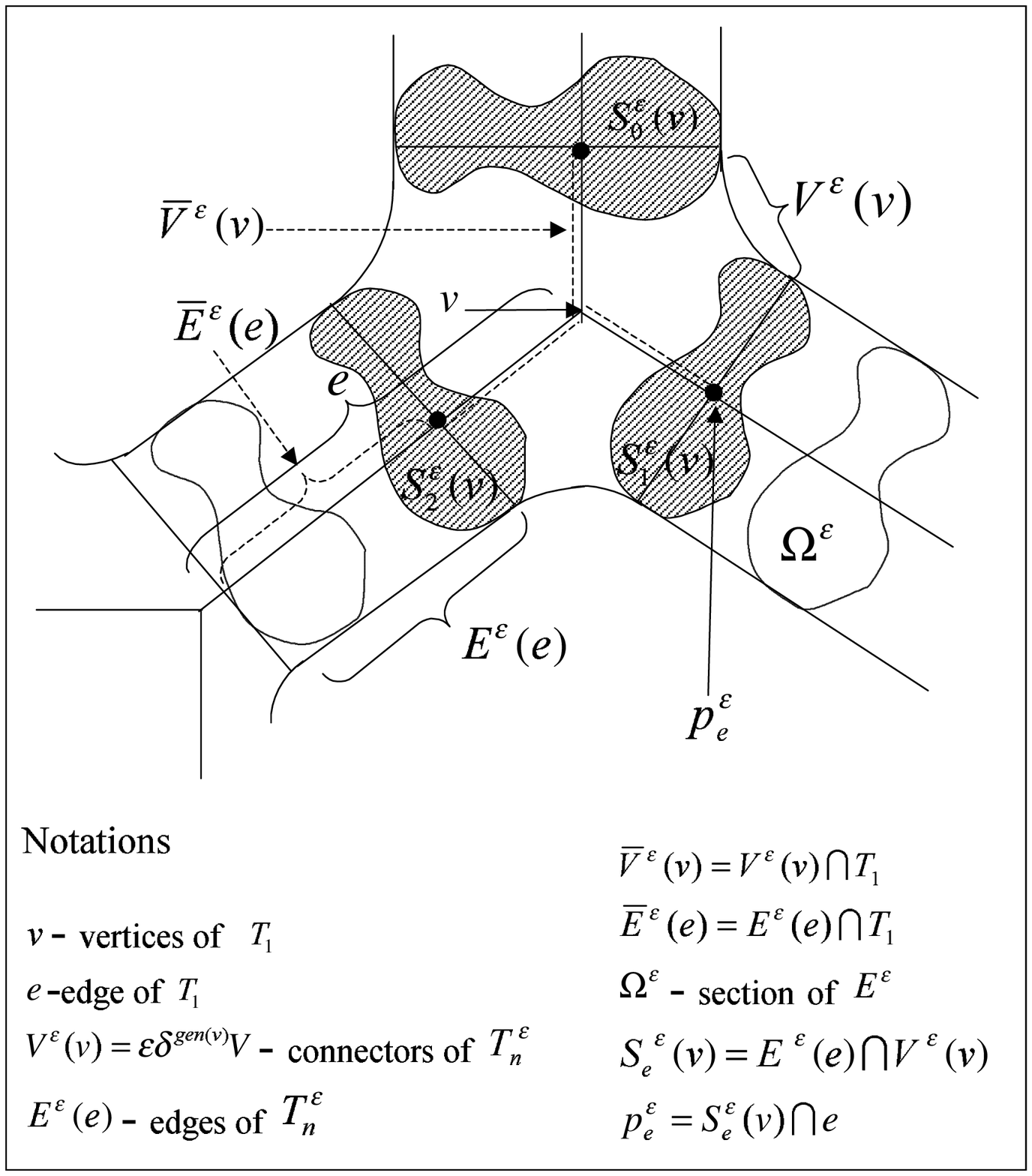}
\caption{Notations of parts of $T_1$ and $T_n^\varepsilon.$}
\label{FigNotationsT2Eps}
\end{figure}

\begin{figure}[t]
\centering
\includegraphics[scale=0.7]{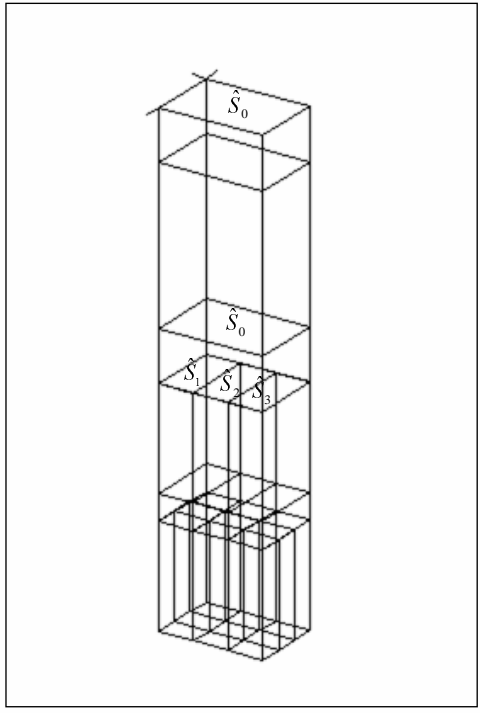}
\caption{The straightened tree, $\hat{T}_3$  for $k=3,N=3$.}
\label{pentagontree3}
\end{figure}

The tree $T_1$ defined above is, in fact, a
{\it combinatorial object},  but we always treat it as a metric
tree or quantum graph. We shall now describe a way to construct an
$N$-dimensional manifold which is an $\eps$-inflation of $T_1$.
For simplicity we shall assume that $T_1$ is radial and regular.
\begin{enumerate}
\item A Lipschitz domain $\Omega\subset \R^{N-1}$ is given.
It corresponds to the (scaled) {\it cross section} of the edges.
We take the origin of $\R^{N-1}$ to be an interior point of
$\Omega$, called the {\it center} of $\Omega$.

\item A Lipschitz domain  $V\subset \R^{N}$ is given. It
corresponds to the (inflated) {\it vertices}. We take the origin
of $\R^{N}$ to be an interior point of $V$, called the {\it
center} of $V$.

 \item $0<\delta <1$ is the scaling factor. The
notation $\delta\Omega$ stands for the scaled domain $\delta\Omega
:= \{ \delta x \mid x\in\Omega\}$. Similarly $\delta V:=\{ \delta
x \mid x\in V\}$. \item The boundary of $V$ contains $k+1$
disjoint sections: One of these sections  is an isometric image of
$\Omega$, denoted by $S_0$. The other $k$ sections are isometric
images of $\delta\Omega$, and denoted by $S_1, \ldots S_{k}$.
\item  The orthogonal projections of
the center of $V$ into $S_0$ and $S_j\subset\partial V$ for $1\leq
j\leq k$   coincide with the isometric image
 of the centers of
$\Omega$ and $\delta\Omega$, respectively.
\end{enumerate}
\noindent
 Next, we define the {\it inflated tree} $T_N^\eps$. For this, let
 us consider a certain embedding of $T_1$ in $\R^{N+1}$.
 We denote this embedding of $T_1$ by the same name,
 $T_1$. It is, in fact, determined by the choice of the inflated
 vertex $V$, to be explained below:
 \begin{enumerate}
 \item[(6)] For each vertex $v$ in the embedded tree $T_1$, the inflated
 vertex is an isometric image of
 $V^\eps(v):=\eps \delta^{\mathrm{gen}(v)}V$  whose center coincides with $v$.
 \item[(7)] Each edge $e\in N(v)$ is
perpendicular to $S_e^\eps(v)$, where $S_e^\eps(v)$ is the
isometric image of the section of $\partial V^\eps(v)$
intersecting the edge $e$.

\item[(8)]   The {\it skeleton} of $V^\eps(v)$ is
 $\overline{V}^\eps(v):= V^\eps(v)\cap T_1$.

 \item[(9)]  For each edge $e$ of the embedded $T_1$, the inflated edge
 is
 $$E^\eps(e)= e\times  S_e^\eps(v)\setminus \cup_v
 V^\eps(v) \ . $$

\item[(10)]    The {\it skeleton} of $E^\eps(e)$ is
 $\overline{E}^\eps(e):= E^\eps(e)\cap T_1$.

 \end{enumerate}
An inflated $2$-dimensional tree is depicted in
Figure~\ref{FigNotationsT2Eps}.
A somewhat degenerate example of an inflated tree is the {\it
straightened tree}, which we denote by $\hat{T}_N$. We use
$\hat{T}_N$ as a canonical representation for $T_N$ in
Section~\ref{ChavelSec}.
\begin{definition}[The straightened tree]\label{stra} {\em
 The  inflated vertex $\hat{V}$ is given by the cylinder
$\hat{\Omega}\times [0,-1]$. The section $\hat{S}_0:=
\hat{\Omega}\times \{0\}$ is the top of $\hat{V}$, and its base
$\hat{\Omega}\times \{-1\}$ consists of $k$ disjoint isometric
copies of $(k)^{-1/N}\hat{\Omega}\times\{-1\}$, corresponding to
the sections $\hat{S}_1, \ldots \hat{S}_k$.  A two-dimensional
straightened tree is depicted in Figure~\ref{pentagontree3}. The
above condition implies that $\hat{\Omega}$ is a box in $\R^{N}$
of a certain type which depend on $k$ and $N$. Indeed, take a box
$\hat{\Omega}$ whose sizes are
$(1,k^{1/N}\!,k^{2/N}\!,\dots,k^{(N-1)/N})$, then $k$ copies of
$(k)^{-1/N}\hat{\Omega}$ exactly cover $\hat{\Omega}$. We may of
course consider also other tilings.}
\end{definition}
\begin{corollary}\label{cortree}
The straightened tree $\hat{T}_N$ can be parameterized by the
cylinder $\hat{\Omega}\times\hat{R}$, where $\hat{R}$ is its
radius (see Figure~\ref{pentagontree3}).
\end{corollary}
%
%
\subsection{Cross sections and functions on $T_1$ and
$T_N^\eps$} \label{SecWeitsOnT1}
There is a natural coordinate system on each of the edges
$E^\eps(e)\subset T_N^\eps$, namely $\vec{x}\in E^\eps(e)\subset
T_N^\eps$ is parameterized as $\vec{x}=\mathbf{x}=(\vec{s},
\theta)$, where $\vec{s}$ is a parameterization of the
corresponding perpendicular section $S_e$ in $\Omega$, scaled by
$\eps\delta^{\mathrm{gen}(e)}$, and $\theta$ is a parameterization
of $e\setminus\cup_v \overline{V}^\eps(v)$. We can also use the
natural parameterization of $V$, scaled by
$\eps\delta^{\mathrm{gen}(v)}$, to describe the coordinate system
in the inflated vertex $V^\eps(v)$. We always take the center of
$V^\eps(v)$ as the origin  $0\in V$.
\begin{enumerate}
\item{We denote by $f_{e}$ the restriction of a function $f$ on
$T_1$ to an edge $e$. In most cases we omit this notation and
write simply $f$ instead of $f_{ e}$.}
\item The function $\rho^*$ is defined on $T_1$ by $\rho^*_{ e}=
\delta^{(N-1)\mathrm{gen}(e)}|\Omega|$.
\item{Let $f$ be a function on $T^1_N$. We denote by
$f^\eps$ the following rescaling of $f$ on $T_N^\eps:$
$$f^\eps(\vec{x})= f_{T_N^\eps}(\vec{x}):=\left\{%
\begin{array}{ll}
    f(\theta,\vec{s}/\eps) & \vec{x}=(\theta,\vec{s})\in
    E^\eps(e),
    \\[2mm]
    f(\vec{x}/\eps) & \vec{x}\in V^\eps(v) ,
\end{array}%
\right.    $$ }
\item The total cross
section of $T_N$  is defined for $t\in T_1\subset T_N$ as
$H(t)=g(t)\rho^*(t)$, where $g$ is the counting function of the
skeleton $T_1$ of $T_N$ and $\rho^*$ as defined in (2) above.
\end{enumerate}

\subsection{Function spaces} \label{secDefA1}
\begin{enumerate}
\item{Let $\rho>0$ be a measurable (weight) function on $T_1$. Denote
$$L_{2,\rho}(T_1)=\left\{f\mid\> f \mbox{ is measurable on } T_1 \mbox{ and
} \int_{T_1} |f|^2\rho\dtet <\infty\right\}.$$ The space
$L_{2,\rho}(T_1)$ equipped with the inner-product $<\!f,g\!>_\rho:
=\int_{T_1} \!\!f\overline{g}\rho\dtet $ is a Hilbert space.}
\item{ $C^1(T_1)$ is the space of continuous functions $f$ on
$T_1$, such that  $f_{ e}\in C^1(e)$ for each edge $e$.} Let
$\rho>0$ be a measurable function on $T_1$. $H^1_{\rho}(T_1)$ is
the completion of the space
$$\left\{f\in C^1(T_1) \mid \sum_{e\in T_1} \int_{e}
\left(|(f_{ e})'|^2+|f_{ e}|^2\right)\rho_{ e}\dtet<\infty
\right\}$$ with respect to the norm
$||f||_{H^1_\rho(T_1)}:=\left[\sum_{e\in T_1} \int_{e} \left(|f_{
e}'|^2+|f_{ e}|^2\right)\rho_{ e}\dtet\right]^{1/2}$.
\item{ $H^1_{0,\rho}(T_1)$ is the completion in $H^1_{\rho}(T_1)$
 of $C_0^1(T_1)$.   For the weight function $\rho^*$, we abbreviate
$H^1_{0,*}(T_1):= H^1_{0,\rho^*}(T_1) $.}
\item{ $H^1(T_N)$ is the completion of the space $$\left\{f\in
C^1(T_N) \mid \int_{T_N}\left(|\nabla
f|^2+|f|^2\right)\mathrm{d}\mathbf{x}<\infty \right\}$$ with
respect to the norm
$||f||_{H^1(T_N)}:=\left[\int_{T_N}\left(|\nabla
f|^2+|f|^2\right)\mathrm{d}\mathbf{x}\right]^{1/2}$.}
\item{ $H^1_0(T_N)$ is the completion in $H^1(T_N)$ of all
functions in $C^1(T_N)$ satisfying $f|_{O\times\Omega_0}=0$. }
\end{enumerate}

\subsection{Laplace and Schr\"odinger operators on $T_1$ and $T_N$} \label{secDefA}

We define a family  of operators on $T_1$ using the standard
definition of operators on $T_1$ (see \cite{Rub1,Sol1}).

Let $W\in L^\infty(T_1)$ be a bounded real valued potential, and
let $\rho_\alpha$ and $\rho_\beta$ be positive bounded
$L^1_{\mathrm{loc}}(T_1)$ weight functions, which satisfy
$\rho_\alpha\asymp \rho_\beta$. In particular,
$H^1_{0,\rho_\alpha}(T_1)$ and $H^1_{0,\rho_\beta}(T_1)$ are
equivalent in the sense that $u\in H^1_{0,\rho_\alpha}(T_1)$ if
and only if $u\in H^1_{0,\rho_\beta}(T_1)$, and there exists a
constant $c>0$ independent of $u$ such that
$$\frac{1}{c}||u||_{H^1_{0,\rho_\alpha}(T_1)}
\leq ||u||_{H^1_{0,\rho_\beta}(T_1)} \leq
c||u||_{H^1_{0,\rho_\alpha}(T_1)}.$$
We denote by
$$E(u,v):=\sum_{e\in T_1} \int_{e}
\left[\frac{\rho_\alpha}{\rho_\beta} (u_{ e})'(\overline{v}_{
e})'+Wu_{ e}\overline{v}_{ e}\right]\rho_\beta \dt $$
the bilinear form on $H^1_{0,\rho_\beta}(T_1)\times
H^1_{0,\rho_\beta}(T_1)$. Without loss of generality, we may
assume that $E\geq 0$ on $C_0^1(T_1)$, so, $E$ is a symmetric and
nonnegative closed bilinear form, and $H^1_{0,\rho_\beta}(T_1)$ is
dense in $L^2_{\rho_\beta}(T_1)$. By Friedrich's extension theorem
(see e.g. Theorem X.23 in \cite{Reed1}) or the First
Representation Theorem (see Theorem VI.2.1 in \cite{Kato}), there
exists a unique selfadjoint operator $L_{\alpha,\beta}$ such that
$\mathrm{Dom}(L_{\alpha,\beta})\subseteq \mathrm{Dom}(E)$ and
$E(u,v)=<L_{\alpha,\beta}u,v>_{\rho_\beta}$ for all $u\in
\mathrm{Dom}(L_{\alpha,\beta})$ and $v\in
H^1_{0,\rho_\beta}(T_1)$. By this theorem, the domain of
$L_{\alpha,\beta}$ is given by:
$$
\mathrm{Dom}(L_{\alpha,\beta})=\{u\in H^1_{0,\rho_\beta}(T_1)\mid
\> |E(u,v)|\leq C|v|_{L^2_{\rho_\beta}(T_1)}\>\> \forall v \in
H^1_{0,\rho_\beta}(T_1) \}
$$
for some constant $C$. Moreover, it is well known (see e.g.
\cite{Rub1}) that the domain of $L_{\alpha,\beta}$ is contained in
the space of all functions $u$ satisfying the following {\em
Kirchhoff conditions}:
\begin{enumerate}
\item{ $u$ is continuous at the vertices (since $
H^1_{0,\rho_\beta}\subset C(T_1)$ ).}
\item{ $\sum_{e\in N(v)} (\rho_\alpha)_{ e} (u_{ e}(v))'=0$ in each
vertex $v\in T_1$.
}
\end{enumerate}

We will call operators of this form {\em width-weighted
operators}, because we will use them for weights $\rho_\alpha$ and
$\rho_\beta$ which are closely related to the width or section
area of $T_N$. Similar operators are also presented by Evans and
Saito in \cite{Eva}.

\begin{remark} \label{RemarkDensityOfDomAalphaBeta} {\em
The domain of the operator $L_{\alpha,\beta}$ is clearly dense in
$H^1_{0,\rho}$ for $\rho=\rho_\alpha$ or $\rho=\rho_\beta$. }
\end{remark}

Finally, the Laplace operator on the tree $T_N$ is defined by the
Friedreich extension of the quadratic form
\begin{equation}\label{En} E_N(u,w):=\int_{T_N}\nabla u\cdot\nabla
\bar{w} \dxx
 \end{equation} for $u,w$ in the space $H^1_0(T_N)$ (see the
definition of $H^1_0(T_N)$  in Section~\ref{secDefA1} \S(5)).
\mysection{Behavior of functions near the vertices}\label{appen1}
 Here we concentrate on a neighborhood of a
vertex (resp. an inflated vertex) in $T_1$ (resp. $T_N^\eps$). For
$T_1$, we shall consider the skeleton $\overline{V}^\eps(v)$
corresponding to a vertex $v$, as defined in
Section~\ref{epsinf}~\S(8). We shall also denote the ``canonical"
skeleton, corresponding to $\eps=1$, by $\overline{V}(v)$.
Occasionally, we shall omit the reference to a particular vertex
$v$ and just denote it as $\overline{V}$.  The end points of
$\overline{V}(v)$ are denoted by $p_e$, where $e\in N(v)$ (see
Figure~\ref{FigNotationsT2Eps}). Recall that $\rho^*$, as defined
in Section~\ref{SecWeitsOnT1} \S(2),  is a positive weight
function on $T_1$, which is constant on each edge.
\par

\begin{enumerate}
\item For each edge $e\in N(v)$ define a nonnegative function $\psi_{(e)}\in
C^1(\overline{V})$ such that $\psi_{(e)}(p_e)=1$ and
$\psi_{(e)}(p_{\tilde{e}})=0$ for $\tilde{e}\not= e$. We also
assume that
\begin{equation}\label{sum1a}\sum_{e\in N(v)} \psi_{(e)} = 1
 \quad \text{ on } \ \overline{V}(v) \  .
\end{equation}
If the skeleton $\overline{V}$ is scaled by $\delta>0$, so
$\overline{V}\rightarrow\delta \overline{V}:= \left\{ \delta
\theta \mid\theta\in\overline{V}\right\}$, where the vertex $v$ is
taken as the origin, then $\psi_{(e)}$ is scaled into
$\psi^{\delta}_{(e)}(x):=\psi_{(e)}( x/\delta)$ for any $x\in
\delta\overline{V}$.
\item Let $V$ be the
``canonical" inflated vertex defined in Section~\ref{epsinf}
\S(2). We choose a family of nonnegative functions $\phi_{(e)}\in
C^1(V)\cap C(\bar{V})$  such that
$$\phi_{(e)}(\mathbf{x})=
  \begin{cases}
    1 &  \theta(\mathbf{x})\in S_e, \\
    0 & \theta (\mathbf{x})\in S_{\tilde{e}}, \mbox{ where } \tilde{e} \neq e,
  \end{cases}
$$
and
\begin{equation}\label{sum1b}\sum_{e\in N(v)} \phi_{(e)}= 1 \qquad
\text{on} \ V   .
\end{equation}
Similarly, if $V$ is scaled by $\delta>0$, so $V\rightarrow\delta
V:= \left\{ \delta x \mid x\in V\right\}$, where the center of $V$
 is taken as the origin, then $\phi_{(e)}$ is scaled into
$\phi^{\delta}_{(e)}(x):=\phi_{(e)}( x/\delta)$ for any $x\in
\delta V$.
\item{ Next, define for each $e\in N(v)$ the quadratic $(k+1)\times( k+1)$
matrices:
$$ {\bf \overline{A}}_{l,m}:= \int_{\overline{V}}(\psi_{(l)})'(\psi_{(m)})'\rho^* \dtet, \ \ \
{\bf A}_{l,m}:=
\int_{V}\nabla\phi_{(l)}\cdot\nabla\phi_{(m)}\dxx,$$ and
$$ {\bf \overline{B}}_{l,m}:=
\int_{\overline{V}}\psi_{(l)}\psi_{(m)}\rho^*\dtet, \ \ \ {\bf
B}_{l,m}:= \int_{V}\phi_{(l)}\phi_{(m)}\dxx.
$$ }
\item{
 Let
$\vec{1}:= \left( 1/\sqrt{k+1}, \ldots 1/\sqrt{k+1}\right)\in
\mathbb{R}^{k+1}$, and for any $\vec{f}\in \mathbb{C}^{k+1}$
denote
\begin{equation} \label{1def} \vec{f} \llcorner\vec{1}:= \vec{f}-
\left(\vec{f}\cdot\vec{1}\right)\vec{1},\end{equation}
where $\cdot$ is the standard inner product in
$\mathbb{C}^{k+1}$.}
\end{enumerate}
%
%

The following Lemma is elementary, but essential for our
 analysis.
\begin{lemma} \label{LemmaAB} The  matrices ${\bf A}$ and ${\bf \overline{A}}$ are nonnegative
definite, and ${\bf B}$ and ${\bf \overline{B}}$ are strictly
positive definite. In particular, there exist constants
$\alpha^{A}>0$, $\alpha^{\overline{A}}>0$, $\alpha^{B}>0$, and
$\alpha^{\overline{B}}>0$, such that
\begin{equation} \label{LemmaABEq1}
\frac{1}{\alpha^{\overline{A}}}|\vec{f} \llcorner\vec{1}|^2
\leq \vec{f}\cdot {\bf \overline{A}}\vec{f}^*
\leq \alpha^{\overline{A}}|\vec{f}
\llcorner\vec{1}|^2,\>\>\>\>\>\>\>\>\>\>
\frac{1}{\alpha^{A}}|\vec{f} \llcorner\vec{1}|^2
\leq \vec{f}\cdot {\bf A}\vec{f}^*
\leq \alpha^{A}|\vec{f} \llcorner\vec{1}|^2,
\end{equation}
 and
\begin{equation} \label{LemmaABEq2}
\frac{1}{\alpha^{\overline{B}}}|\vec{f}|^2
\leq \vec{f}\cdot {\bf \overline{B}}\vec{f}^* \leq
 \alpha^{\overline{B}}|\vec{f}|^2,\>\>\>\>\>\>\>\>\>\>
 \frac{1}{\alpha^{B}}|\vec{f}|^2
\leq \vec{f}\cdot {\bf B}\vec{f}^* \leq
 \alpha^{B}|\vec{f}|^2
 \end{equation}
for all $\vec{f}\in \mathbb{C}^{k+1}$, where $\vec{f}^*$ denotes
the complex conjugate of $\vec{f}^t$.
\end{lemma}
\begin{proof} The non-negativity (resp. positivity) of ${\bf
A}$ and  $\overline{\bf A}$  (resp. ${\bf B}$ and  $\overline{\bf
B}$) follows from the corresponding definitions, while
(\ref{LemmaABEq1}) follows from (\ref{sum1a}) and  (\ref{sum1b}).
\end{proof}
Let us introduce the following functionals on
$H^1(\overline{V})$:
\begin{equation}\label{LemmaDirichletSolExistOnVBarDef1}
    \overline{I}^{\overline{V}}_\gamma[g]:=\int_{\overline{V}}(|g'|^2+\gamma
    |g|^2)\rho^*\dtet
    \qquad\mbox{ for }\>\>\gamma=0,1,
\end{equation}
and for $\vec{f}\in \mathbb{C}^{k+1}$ let us denote:
\begin{equation}\label{LemmaDirichletSolExistOnVBarDef2}
    \overline{\mathcal{A}}_{\overline{V},\vec{f}}=\{g\in
    H^1(\overline{V}) \mid g(p_e)=f_e, \ \ e\in N(v) \ \}.
\end{equation}
\begin{lemma}\label{LemmaDirichletSolExistOnVBar}
Using the notations  \eqref{LemmaDirichletSolExistOnVBarDef1} and
\eqref{LemmaDirichletSolExistOnVBarDef2}, we have for $\gamma=0,1$
that
$$\overline{J}^{\overline{V}}_\gamma[{\vec{f}}]
:=\inf_{g\in\overline{\mathcal{A}}_{\overline{V},{\vec{f}}}}\overline{I}^{\overline{V}}_\gamma[g]
$$
 is attained by a unique function $h$, which
solves the Dirichlet problem
\begin{equation} \label{LemmaDirichletSolExistOnVBarEq1}
    -h''+\gamma h=0 \quad \mbox{ in } \overline{V}\cap e
    ,\>\>\>h(p_e)=f_e
    \quad\forall e\in N(v),
\end{equation}
and satisfies Kirchhoff's conditions
\begin{equation}\label{Kirk} \sum_{e\in N(v)} \rho^*_{ e}
h_e'(v) = 0. \end{equation}
\end{lemma}
\begin{proof} The existence of minimizers $u$ for
$\overline{I}^{\overline{V}}_0$ and
$\overline{I}^{\overline{V}}_1$, which satisfy
\eqref{LemmaDirichletSolExistOnVBarEq1} is standard (see e.g. the
proof in \cite[Theorem 2, pp. 448--449]{Eva}).

We need to prove that the minimizer $u$ of
$\overline{I}^{\overline{V}}_\gamma$ satisfies Kirchhoff's
derivatives condition. To this end, let $v\in C_0^1(\overline{V})$
and $0\neq\epsilon\in \mathbb{R}$. Since $u$ is a minimizer,
$I_\gamma[u] \leq I_\gamma[u+\epsilon w], $ and therefore,
$$\int_{\overline{V}}(u'w'+\gamma uw)\rho^*\dtet=0. $$
By elliptic regularity $u\in C^2(V\cap e);$ Moreover, $u$ is
continuous in $\overline{V}$. Recall that  $\rho^*$ is constant on
each edge, therefore, $ -u^{''}+\gamma u=0$ on $\overline{V}\cap
e$. Thus, \begin{multline} 0=\sum_{e\in
N(v)}\int_{\overline{V}\cap e}(u'w'+\gamma uw)\rho^*\dtet\\
=\left.\sum_{e\in N(v)}\rho^*_{ e} (u_{ e})'w_{ e}
\right|_{p_e}^{v_e}+ \sum_{e\in N(v)} \int_{\overline{V}\cap
e}(-u''+\gamma u)w\rho^*\dtet   = w(v)\sum_{e\in N(v)}\rho^*_{ e}
 u_{e}'(v) .
\end{multline} The uniqueness of the minimizers of
$\overline{I}^{\overline{V}}_0$ and
$\overline{I}^{\overline{V}}_1$ follows since both are minima of
strictly convex functionals on the underlying domains.
\end{proof}
\begin{lemma}\label{LemmaFormOnVBarj}
There exist $\beta^{\overline{A}}>0$ and $\beta^{\overline{B}}>0$
such that for all $\delta>0$
\begin{equation}\label{LemmaFormOnVBarjEq1}
    \overline{I}^{\delta\overline{V}}_0[\vec{f}]\geq
    \delta^{-1}\beta^{\overline{A}}|\vec{f} \llcorner\vec{1}
    |^2,
\end{equation}
 and
\begin{equation}\label{LemmaFormOnVBarjEq2}
    \overline{I}_1^{\delta\overline{V}}[\vec{f}]\geq
    \delta^{-1}\beta^{\overline{B}}
(|\vec{f} \llcorner\vec{1}
    |^2 + \delta^2 |\vec{f}|^2).
\end{equation}
\end{lemma}
\begin{proof} In the following, we use the notations
introduced in Lemma \ref{LemmaDirichletSolExistOnVBar}, and in
\eqref{LemmaDirichletSolExistOnVBarDef1} and
\eqref{LemmaDirichletSolExistOnVBarDef2}. Consider the case
$\delta=1$ first. Let $\{\vec{\sigma}_e\}$ be the standard basis
vectors  in $ \mathbb{C}^{k+1}$, where $e\in N(v)$. Let
$h_{(e)}\in H^1(\overline{V})$  be  the unique minimizer of
$\overline{J}^{\overline{V}}_\gamma[{\sigma_e}]$. By
Lemma~\ref{LemmaDirichletSolExistOnVBar} it follows that
$$ \overline{J}^{\overline{V}}_\gamma [ \vec{f}] =
\overline{I}^{\overline{V}}_\gamma \left[ \sum_{e\in N(v)}
f_eh_{(e)}\right] =\sum_{e,\tilde{e}\in N(v)} f_e
f_{\tilde{e}}\int_{\overline{V}}
 \left[h_{(e)} h_{(\tilde{e})}' + \gamma h_{(e)}h_{(\tilde{e})}\right]\rho^* \dtet, $$
where each $h_{(e)}$ satisfies
$$-h_{(e)}''+ \gamma h_{(e)}=0 \quad \mbox{ in } \overline{V},\quad h_{(e)}(p_e)=1,
\ \ h_{(e)}(p_{\tilde{e}})=0
 \>\>\>\forall
\tilde{e}\not= e .$$   Let $\gamma=0$. By Lemma
\ref{LemmaDirichletSolExistOnVBar}, $\overline{J}_0[\vec{f}]$ is
attained uniquely by the harmonic function $h$ which solves the
corresponding Dirichlet problem (and satisfies Kirchhoff's
conditions). In particular, it depends only on $\vec{f}$ and the
domain $\overline{V}$. Since each solution $h$ satisfying
$h(p_e)=f_e$ can be presented uniquely by $h=\sum_{e\in
N(v)}f_eh_{(e)}$, it follows that $\overline{J}_0[f]$ is a
bilinear form. Clearly, it is a nonnegative $k+1$ dimensional form
whose kernel contains only constant multiplicities of ${\vec{1}}$
for which the unique solution of the Dirichlet problem is
constant. Therefore, it is equivalent to all nonnegative forms
with such a kernel, and in particular, to $|\vec{f}
\llcorner\vec{1}|^2$.

The proof for the case $\gamma=1$  is similar except for replacing
the Laplace operator by the operator
$-\mathrm{d}^2/\mathrm{d}\theta^2+1$ and  $|\vec{f}
\llcorner\vec{1}|^2$ by $|\vec{f}|^2$.

\par Now, if $\delta<1$ and
$\gamma=0$ we observe that the harmonic minimizers $h_{(e)}$ are
scaled into $h_{(e)}(\cdot/\delta)$, and
$$ \int_{\delta\overline{V}} h'_{(e)}(x/\delta)
h'_{(\tilde{e})}(x/\delta)\rho^*\dtet=
\delta^{-1}\int_{\overline{V}} h_{(e)}' h_{(\tilde{e})}'\rho^*
\dtet .
$$ For $\delta <1$ and $\gamma=1$, we use similar
scaling argument to obtain (\ref{LemmaFormOnVBarjEq2}).
\end{proof}

We wish to prove now the analog of Lemma \ref{LemmaFormOnVBarj}
for the $N$-dimensional case. Consider the following functionals
for $\gamma=0,1:$
\begin{equation}\label{LemmaDirichletSolExistOnVjDef1}
    I_\gamma [g]:=\int_{V}(|\nabla g|^2+\gamma
    |g|^2)\dxx.
\end{equation}
For all $h\in H^1(V)$ and $0\leq j\leq k$ we denote the {\it
average} of $h$ on the section  $S_j\subset \partial V$ by
\begin{equation}\label{LemmaDirichletSolExistOnVjDef2}
     P_j (h) :=\frac{1}{|S_j|}\int_{S_j} h \dss
\end{equation}
(see Section~\ref{epsinf} \S(4)). For $\vec{F}\in
\mathbb{C}^{k+1}$ we define
\begin{equation}\label{LemmaDirichletSolExistOnVjDef3}
     \mathcal{A}_{V,\vec{F}}:=\left\{g\in H^1(V)\mid P_j(g)=F_j\quad \forall
    j=0,...,k\right\}.
\end{equation}
\begin{lemma} \label{LemmaDirichletSolExistOnVj}
Let $\vec{F}\in \mathbb{C}^{k+1}$. Using the above notations, we
have for $\gamma=0,1$ that
$$J_\gamma[\vec{F}]:=\inf_{ g\in\mathcal{A}_{V,\vec{F}}}I_\gamma[g] $$
is attained by a function $h$, which is the unique solution of the
problem
\begin{equation} \label{LemmaDirichletSolExistOnVEq1}
    -\Delta h+\gamma h=0 \>\>\mbox{ in }\>\>V,\>\>\>\>\>h\in
    \mathcal{A}_{V,\vec{F}}\, ,
\end{equation}
and satisfies weakly the mixed boundary conditions
\begin{equation} \label{LemmaDirichletSolExistOnVEq1_5}
    \frac{\partial h}{\partial n}=0 \>\>\mbox{ on }\>\>\partial
    V\backslash\cup_{j=0}^{k}S_{j}, \>\>\>\>\>\mbox{ and }\>\>\>\>\>
    \frac{\partial h}{\partial n}=\kappa_j \>\>\mbox{ on }\>\>S_{j},
\end{equation}
where $\kappa_{j}$ for $j=0,...,k$ are uniquely determined
constants.
\end{lemma}

\begin{proof} The proof of
\eqref{LemmaDirichletSolExistOnVEq1} for the case $\gamma=1$ is
standard.  Indeed,  let $\{w_i\}_{1=1}^\infty$ be a minimizing
sequence satisfying
$\lim_{i\rightarrow \infty}I_1[w_i]=J_1[\vec{F}]$. Then
$\{w_i\}_{i=1}^\infty$ is bounded in $H^1(V)$.  Therefore, there
exists a subsequence $\{w_i\}$ and a function $v\in H_1(V)$ such
that $w_i\rightharpoonup v \>\>\mbox{ in }\>\>H^1(V)$.

Since $P_j(f)$ is a continuous functional on $H^1(V)$ in the
strong topology, it is also continuous in the weak topology. In
particular, $\mathcal{A}_{V,\vec{F}}$ is closed in the weak
topology of $H_1(V)$ so $v\in \mathcal{A}_{V,\vec{F}}$ . The lower
semicontinuity of $I_1$ implies that $I_1[v]=J_1(\vec{F})$.
Moreover, $v$ is unique because $I_1$ is convex.

It remains to prove that $v$ satisfies the boundary conditions in
\eqref{LemmaDirichletSolExistOnVEq1_5}.  Since $v$ is a minimizer
in $\mathcal{A}_{V,\vec{F}}\,$, it follows that
$$0=\int_{V}(\nabla v\cdot\nabla
w+vw)\dxx =\int_{\partial V\backslash
\cup_{j=0}^{k}S_{j}}w\frac{\partial v}{\partial n}\,\mathrm{d}\xi
+\int_{ \cup_{j=0}^{k}S_{j}}w\frac{\partial v}{\partial n} \dss .
$$
for any $w\in \mathcal{A}_{V,\overrightarrow{0}}$.  The first term
in the last expression is thus  zero only  if $\partial v/\partial
n=0$ on $\partial V\backslash \cup_{j=0}^{k}S_{j}$ in the weak
sense. Since the average of the test function $w$ is zero on each
sector ($P_j(w)=0$), the second term is zero
 if ${\partial v}/{\partial n}=\kappa_{j}$ (in the weak sense) on $S_{j}$.
Finally, the multipliers $\kappa_{j}$ are uniquely determined due
to the uniqueness of $v$ for any $\vec{F}$.

The proof of (\ref{LemmaDirichletSolExistOnVEq1_5}) for the case
$\gamma=0$ is similar, except that we have to prove the bound in
$L^2(V)$ of the minimizing sequence.  Since $V$ is a bounded
Lipschitz domain, by \cite[Theorem 5.5.1, and the remark in p.
286]{Maz}, the embedding $H^1(V)\rightarrow L^2(V)$ is compact.
Hence,  the spectrum of Helmholtz operator with the Neumann
boundary condition for such domains is discrete. Its first
eigenvalue is $1$, and is a simple eigenvalue corresponding to the
constant ground state. Hence, the Poincar\'{e} inequality
\begin{equation}\label{LemmaDirichletSolExistOnVBarEq4}
    \int_{V}|v|^2\mathrm{d}\mathbf{x}\leq \Lambda_2^{-1} \int_{V}|\nabla v|^2\mathrm{d}\mathbf{x}.
\end{equation}
holds for all functions $v$ perpendicular to the constant in
$H_1(V)$, where $\Lambda_2$ is the second eigenvalue of the
Neumann Laplacian on $V$.

We now repeat the argument for the case $\gamma=1$, but restrict
our domain to the domain of all functions in
$v\in\mathcal{A}_{V,\vec{F}}$ which are perpendicular to the
constant. The minimizer $u$ obtained in this way satisfies $P_j
(u) = F_j+\kappa$ for some $\kappa\in \R$ and $j\in \{0,\ldots,
k\}$. Then $u-\kappa\in \mathcal{A}_{V,\vec{F}}$.
\end{proof}

Let now $\delta>0$, and set $\delta V:=\{ \delta x \ ; \ x\in V\}$
the scaled inflated vertex, where we assume (as usual) that the
center of $V$ is in the origin. The sections of $\delta V$ are
scaled accordingly, and we denote them by $\delta S_j$, $0\leq
j\leq k$. We define, correspondingly, the averaging operator on
$\delta S_j$ for $h\in H^1(\delta V)$:
\begin{equation}\label{deltaP}
     P^\delta_j (h):=\frac{1}{\delta^{N-1}|S_j|}\int_{\delta S_j} h \dss ,
\end{equation}
and
\begin{equation}\label{Adelta1}
     \mathcal{A}^{\delta}_{V,\vec{F}}:=
     \left\{g\in H^1(\delta V)\mid P^\delta _j(g)=F_j\quad\forall
    j=0,...,k\right\}
   \ .
\end{equation}
Using Lemma~\ref{LemmaDirichletSolExistOnVj}, the following lemma
is proved analogously  to the proof of the second part of
Lemma~\ref{LemmaFormOnVBarj}.
\begin{lemma}\label{LemmaFormOnVj}
 There exist $\beta^A>0$ and $\beta^B>0$ such that for all $f\in
\mathcal{A}^{\delta}_{V,\vec{F}} $
\begin{equation}\label{LemmaFormOnVjEq1}
    \int_{\delta V}|\nabla f |^2 \dxx\geq  \delta^{N-2}\beta^A|\vec{F} \llcorner\vec{1}
|^2, \end{equation} and
\begin{equation}\label{LemmaFormOnVjEq2}
    \int_{\delta V} (|\nabla f|^2 + |f|^2)\dxx \geq  \beta^B\left(\delta^{N-2}|\vec{F} \llcorner\vec{1}
|^2 + \delta^N|\vec{F}|^2\right)  .
\end{equation}
\end{lemma}

\mysection{Discreteness  of the spectrum on $T_1$ and $T_N$}
\label{SectionBackground}

In this section we study the discreteness of the spectrum of
width-weighted operators on $T_1$ and Schr\"{o}\-din\-ger
operators on $T_N$.

\subsection{Discreteness of the spectrum for weighted operators on $T_1$}
In \cite{Car}, Carslon has shown that the spectrum of the
Laplacian on a connected metric graph $G$ of finite volume which
has a compact completion $\overline{G}$ is purely discrete.
Solomyak \cite{Sol1} has extended Carlson's result to regular
trees of  a finite radius $R$.
\begin{theorem}[Solomyak \cite{Sol1}] \label{SolTh1}
    Let $T_1$ be a radial tree such that $R(T_1)<\infty$ and   its branching function  is uniformly bounded.
        Let $W(x)$ be a radially symmetric measurable real
    valued function which is bounded below. Then the spectrum of
    $-\Delta+W$ on $T_1$ is purely discrete.
\end{theorem}

\begin{proof}[Outline of Solomyak's proof] Solomyak constructed a family of weighted operators
    $\{A_{W,v}\}$ which are defined on the intervals $[t_v,R)\subseteq
    \mathbb{R}$, where $t_v$ is the distance of a vertex $v$ from
    the root $O$.
    The operators $A_{W,v}$ are defined as the selfadjoint operators in
    $L_{g}^2(t_v,R)$, associated with the quadratic form
    \begin{equation} \label{aVk}
        a_{W,v}[u]:=\int_{t_v}^{R}\left[|u'(t)|^2+W(t)|u(t)|^2\right ]g(t)\dt\qquad u\in C^\infty_0(t_v,R),
            \end{equation}
            where $g$ is the counting function. Using a
decomposition of functions in $H(T_1)$ into symmetric functions on
subtrees \cite{Naim1} (which implies the spectral decomposition of
the Laplacian to these operators), Solomyak showed the equivalence
between the discreteness of the
    spectrum of the Laplacian on $T_1$ and the discreteness of the
    spectrum of $A_{W,v}$ on $[t_v,R)$ for all vertices $v\in
    T_1$. Using a theorem of Birman and Borzov \cite{Bir} and a certain
change of variables, it is then shown that all the operators
$A_{W,v}$ have discrete spectra. The proof of this part relies on
the monotonicity of the counting function $g$.
\end{proof}

The basic ingredient in Solomyak's proof, namely the spectral
decomposition  into the subspaces of functions which are symmetric
on subtrees, still holds if one adds weight functions which are
symmetric in generations (see \cite{Naim1,Sol1} for details). The
Schr\"odinger-type operators we consider in this section are
defined on the {\it weighted} tree $T_1$  and involve a pair of
symmetric weight functions $\rho_\alpha, \rho_\beta$ and a
symmetric potential $W$:
\begin{equation}\label{Lab}\overline{L}_{\alpha,\beta}u:=
-\rho_\beta^{-1} \frac{\mathrm{d}}{\mathrm{d}t} \left(\rho_\alpha
\frac{\mathrm{d}u}{\mathrm{d}t}\right) + Wu.\end{equation} The
spectral decomposition of $\overline{L}_{\alpha,\beta}$ is
obtained by reducing these operators to the space of functions
which are symmetric on all subtrees. The restriction of
$\overline{L}_{\alpha,\beta}$
 to the symmetric subtree
$T_{1,v}$ with a root $v\in T_1$ are obtained by the quadratic
form
\begin{equation} \label{OuraVk}
    a_{\alpha,\beta,W,v}[u]=\int_{t_v}^{R}
    \left[|u'(t)|^2\rho_\alpha+W(t)|u(t)|^2\rho_\beta\right]g(t)\dt.
\end{equation}
and the associated operator in $L_2([t_v, R))$ is denoted by $
A_{\alpha,\beta,W,v}$. To extend the result of Solomyak to the
weighted tree we should show that $ A_{\alpha,\beta,W,v}$  has a
discrete spectrum for each vertex $v\in T_1$. Even though
(\ref{OuraVk}) seems very close to (\ref{aVk}), the counting
function $g$ in (\ref{aVk}) is replaced by $g\rho_\alpha$ and $g
\rho_\beta$ in (\ref{OuraVk}), and these functions are not
necessarily monotone. We prove the discreteness of
$\overline{L}_{\alpha,\beta}$ under the weaker condition that
$g\rho_\alpha$ and $g \rho_\beta$ are uniformly bounded from
below:
\begin{theorem}\label{LemmaT1DiscLew}
    Let $T_1$ be a one-dimensional tree, whose radius $R$ is
    finite. Assume that  $0<\rho<1$ is a symmetric weight function on $T_1$,
that $\rho_\alpha\asymp\rho$ and $\rho_\beta\asymp\rho$ are
symmetric weight functions. Suppose that there exists a constant
$C>0$ so that
    \begin{equation}\label{LemmaT1DiscLewEq1}
        Cg(s)\rho(s)<g(t)\rho(t) \qquad \mbox{ for all } \ \ \
       s\leq t \leq R(T_1)  .
    \end{equation}
    Then the spectrum of the width-weighted operator $\overline{L}_{\alpha,\beta}$ on
    $T_1$ is purely discrete.
\end{theorem}
We use the following general Lemma of Lewis.
\begin{lemma}[{\cite[Lemma 1]{Lew}}] \label{LewLemma}
    Let $D$ be a domain in $\mathbb{R}^N$.
    Let $\mathrm{h}$ be a strictly positive symmetric closed
    form whose domain $H_{\mathrm{h}}(D)$ is dense in the Hilbert space
    $L_w^2(D)$ for a positive weight function $w$ on $D$.\\
    Suppose that $D$ is the union of an increasing sequence
    of open sets $\{D_j\}$, for which the identity injection
    $i_j:H_{\mathrm{h}}(D_j)\rightarrow L^2_w(D_j)$ is compact.\\
    If there is a positive-valued function $p(x)$ on $D$ and a
    sequence of positive numbers $\eps_j\rightarrow 0 $ as
    $j\rightarrow \infty$ such that $$w(x)p(x)^{-1}<\eps_j \qquad\mbox{ for
    almost every } x\in D \backslash D_j,$$ and
    \begin{equation} \label{LewEq}
        \int_{D \backslash D_j} p(x)|u(x)|^2\mathrm{d}x\leq
        \mathrm{h}[u,u] \qquad \mbox{ for all } u\in H_{\mathrm{h}},
    \end{equation} then
    the selfadjoint operator on $L_w^2(D)$ associated with the
    Friedrich extension of $\mathrm{h}$ has a purely discrete
    spectrum.
\end{lemma}

\begin{proof}[Proof of Theorem \ref{LemmaT1DiscLew}] We only need
to  show that for any $v\in T_1$ the operator
$A_{\alpha,\beta,W,v}$ associated with (\ref{OuraVk}) has  a
discrete spectrum. Evidently, it is enough to show it for $v=O$.
For this, we use Lemma~\ref{LewLemma} with the quadratic form
$h=a_{\alpha, \beta}$ on $L^2_{\rho_\beta}=L^2((0,R),\rho_\beta
\dt)$. We set $D=[0,R)$, $D_j=[0, t_j)$. We denote
$$p(\theta):=\frac{\rho(\theta)g(\theta)}{ R(R-|\theta|)} \ . $$
Since $0<\rho<1$ and $g\rho$ satisfies  \eqref{LemmaT1DiscLewEq1},
it follows that  $p$ satisfies the assumptions of
Lemma~\ref{LewLemma}. By our assumptions $\rho_\alpha\asymp
\rho_\beta\asymp \rho$, therefore it is sufficient to prove that
for all $u\in C^1_0([0,R))$ and $0<j<R$ we have
$$
\int_j^R p(\theta)|u(\theta)|^2\,\mathrm{d}\theta
\leq
C\int_0^R|u'(\theta)|^2\rho(\theta)g(\theta)\,\mathrm{d}\theta.
$$
In fact, for any $j<\theta<R$
$$|u(\theta)|^2
=\left| \int_{\theta}^{R}u'(t)\mathrm{d}t\right|^2\leq
|R-\theta|\int_\theta^R |u'(t)|^2\mathrm{d}t.
$$
Then \begin{multline*} \int_j^R p(\theta)|u(\theta)|^2
\,\mathrm{d}\theta \leq \int_j^R |R-\theta|p(\theta)
\left[\int_\theta^R|u'(\zeta)|^2\,\mathrm{d}\zeta\right]\,\mathrm{d}\theta
  \\[2mm] \leq \frac{1}{R}\int_j^R \left[\rho(\theta)g(\theta)\int_\theta^R
  |u'|^2\,\mathrm{d}\zeta\right]\,\mathrm{d}\theta
 \leq \frac{C}{R}\int_j^R
 \left[\int_\theta^R \rho g|u'|^2\,\mathrm{d}\zeta \right]\,\mathrm{d}\theta
 \\[2mm]  \leq C\int_0^R \rho g |u'|^2
\,\mathrm{d}\theta
  \ . \end{multline*}
  \end{proof}


\subsection{Discreteness of the spectrum for operators on
$T_N$} \label{ChavelSec} As we have mentioned, we are interested
in spectral properties of Schr\"{o}dinger operators on the
$N$-dimensional tree $T_N^\eps$. It is well known that the
Laplacian on a compact manifold with a smooth boundary, and with
standard (regular) boundary conditions has a pure point spectrum.
However, since we wish to address also the problem of the
discreteness of the spectrum for nonsmooth trees with an infinite
volume, we cannot implement the classical theory. Instead, we use
Lemma~\ref{LewLemma} to prove the discreteness of the spectrum of
Schr\"{o}dinger operator on $T_N^\eps$ with a finite radius.

Recall Definition~\ref{stra} of the straightened  tree
$\hat{T}_N$. By Corollary~\ref{cortree} we can assign $\hat{T}_N$
a global coordinate system to the tree, namely $(\vec{s},\theta)$,
where $\vec{s}\in \hat{\Omega}$ and $\theta\in [0, \hat{R})$.
 We pose the following assumption.
\begin{assumption}\label{asstree0}
{\em There exists a $C^1$-diffeomorphism
$\mathcal{G}:{T}_N\rightarrow \hat{T}_N$. We denote by
$\mathcal{F}$ its inverse, so that,
$\mathcal{F}:\hat{T}_N\rightarrow T_N$. Denote by $J$ the Jacobian
of $\mathcal{F}$. We assume that there is a constant $C>0$ such
that \bea\label{assum1} \left|\frac{\partial
\mathcal{F}(\vec{s},\theta)}{\partial \theta}\right|\leq C \qquad
\forall (\vec{s},\theta) \in \hat{T}_N,\\[3mm]
 \label{assum2} 0<J(\vec{s},\theta_1)\leq C J(\vec{s},\theta_2) \qquad \forall
\theta_1\leq \theta_2 .\eea
 } \end{assumption} We have in mind the following two-dimensional
example.
\begin{example} {\em \label{Example2DTreedisc}
Let $T_2$ be a two-dimensional binary symmetric tree constructed
by gluing rectangles and triangles (see
Figure~\ref{pentagontrans}). Let the length of a rectangle in
generation $j$ be $r^j$ and its width be $d^j$, where
 $r,d\in(0,1)$. Clearly, $T_2$ can be embedded in $\mathbb{R}^3$ to avoid
overlapping of the edges. Notice that such a tree may have an
infinite area (though its radius is finite). Indeed,  the area of
such a tree is given by $\sum_{j=1}^\infty [(2dr)^j+\beta d^{2j}]$
for some constant $\beta$, hence for any choice of $r<1,d<1$ such
that $rd>1/2$, the area is infinite. Let us denote by $P_j$ the
pentagon constructed by gluing a rectangle and triangle in the $j$
generation, and by $P_{j,l}$ for $l=1,2$ its partition into two
symmetric quadrangles. We assume that the coordinates of the
vertices of the quadrangle $P_{j,l}$, $(x_{1,a},x_{2,a})$ for
$a=1,...,4$, are given (up to translations) by
$(0,0),(d^j,0),(d^j,r^j),(0,r^j+cd^j)$ respectively for a constant
$c$.

Let $p={(R-1)}/{R}$, where $R$ is the radius of the original tree
($p$ is chosen such that $\sum_{j=0}^\infty p^j=R$). In
particular, $r<p<1$.
A transformation of a rectangle whose vertices are at
$(0,0),(1/2^j,0),(1/2^j,p^j),(0,p^j)$  onto $P_{j,l}$ can be
written in the form
$$
x_1(\theta,s)=(2d)^js, \>\>\>\mbox{ and }\>\>\>
x_2(\theta,s)=\frac{r^j}{p^j}\theta+
c\frac{d^j}{p^j}\theta-c\frac{2^jd^j}{p^j}s\theta.
$$
An elementary  calculation shows that
if $d\leq p$, then $ \left|{\partial x_2}/{\partial \theta}\right|
\leq 1+c $. Note that ${\partial x_1}/{\partial s}=(2d)^j$ is not
bounded for $d>1/2$, which means that the total width of the tree
is unbounded. However, the condition $d>1/2$ ensures the
possibility of gluing together the connectors and the edges of
this tree.
\begin{figure}[t]
\centering
\includegraphics[scale=0.7]{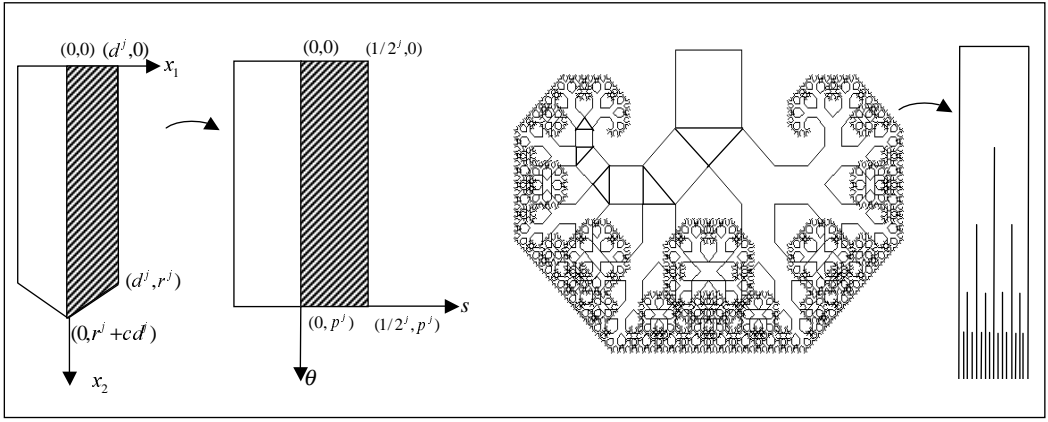}
\caption{The transformation of $T_2$ to $\widehat{T}_2$.}
\label{pentagontrans}
\end{figure}
}
\end{example}
 \begin{assumption}\label{asstree} {\em
Let $\hat{V}\subset\hat{\Omega}\times [0,1]$ be an inflated vertex
of the straightened tree, where $\hat{S}_0:=\hat{\Omega}\times
\{0\}$, $\hat{S}_j\cong k^{-1/N}\hat{\Omega}\times\{1\}$, $1\leq
j\leq k$ the corresponding sections. Let $V$ be the inflated
vertex of a given tree $T_N$, and $S_j\subset\partial V$, $0\leq
j\leq k$ the corresponding sections. We assume that there exists a
$C^1$-diffeomorphism
$\mathcal{F}=\mathcal{F}(\vec{s},\theta):\hat{V}\rightarrow V$ so
that $\mathcal{F}(\hat{S}_j)=S_j$ for $0\leq j\leq k$, and such
that
$$\left| \frac{\partial \mathcal{F}}{\partial \theta}\right|<C,
\quad \mbox{ and } \quad
 0<J(\vec{s},\theta_1)\leq C J(\vec{s},\theta_2) \qquad \mbox { if }
\theta_1\leq \theta_2   $$ hold on $V$ for some $C>0$, where $J$
is the Jacobian of $\mathcal{F}$. }
\end{assumption}
\begin{remark}\label{lemassump}
Assumptions~\ref{asstree} imply Assumptions~\ref{asstree0}.
\end{remark}
\begin{theorem} \label{lemmaDiscreteSpecT2}
Under Assumptions~\ref{asstree0} (resp.
Assumptions~\ref{asstree}), the Laplace operator on $T_N$ as
defined in Section~\ref{secDefA}, has a purely discrete spectrum.
\end{theorem}
{\bf Proof.}  Let  $\mathcal{G}:T_N\rightarrow \hat{T}_N$ be the
inverse $C^1$-mapping of $\mathcal{F}$ which is defined in
Assumptions~\ref{asstree0}, and set $\mathcal{G}(\mathbf{x}):=
(\theta(\mathbf{x}), \vec{s}(\mathbf{x}))\in \hat{T}_N$. Denote by
$J$ the Jacobian of $\mathcal{F}$. Let $\hat{T}_{N,j}\subset T_N$
be the finite subtree
$$ \hat{T}_{N,j}:=\left\{ (\theta,\vec{s})\in \hat{T}_N \mid
\theta < \theta_j \right\},$$ where $\theta_j \nearrow \hat{R}$,
and $\hat{R}$ is the radius of $\hat{T}_N$. Let
\begin{equation}\label{Tnj} T_{N,j}:= \mathcal{F}
\left(\hat{T}_{N,j}\right),
\end{equation}
and set
$$p(\mathbf{\mathbf{x}}):=\frac{1}{C^2\hat{R}|\hat{R}-\theta({\mathbf{x}})|} \ . $$
We wish to use Lemma~\ref{LewLemma} with $D_j\equiv T_{N,j}$. This
Lemma requires the compactness of the identity injection
$i_k:H^1(D_j)\rightarrow L^2(D_j)$. Although the boundary of
$D_j=T_{N,j}$ is not $C^1$, this injection is still compact.
Indeed, the embedding $i:H^1(D)\rightarrow L^2(D)$ is compact for
a bounded domain $D$ which has the (inner) cone property (see
\cite[Theorem 5.5.1]{Maz}, and the remark on p. 286 therein).

By Lemma~\ref{LewLemma}, it is sufficient to prove for the
Laplacian that
\begin{equation}\label{Eq1DiscreteTnProof}
\int_{T\backslash
T_{N,j}}p(\mathbf{x})|u(\mathbf{x})|^2\mathrm{d}\mathbf{x}
\leq \int_{T_N}|\nabla u|^2\mathrm{d}\mathbf{x}
\end{equation}
for all $u\in C^1(T_N)$ that vanish on the `top' of $T_N$ and
outside $ T_{N,j}$ for some $j\geq 1$.  Let $u$ be such a test
function,  and let $v(\theta,\mathbf{s})=u(\mathbf{x})$. Then
\begin{equation}\label{Eq2DiscreteTnProof}
|u(\mathbf{x})|^2\!\!=\!\!|v(\theta,\mathbf{s})|^2
=\left|\int_{\theta}^{\hat{R}} \frac{\partial v}{\partial
\vartheta}\mathrm{d}\vartheta\right|^2
\leq |\hat{R}-\theta|\int_{\theta}^{\hat{R}} \left| \frac{\partial
v}{\partial \vartheta}\right|^2 \mathrm{d}\vartheta.
\end{equation}

 Using the definition
of the function $p$, \eqref{assum1}, \eqref{assum2},
\eqref{Eq2DiscreteTnProof}, and Fubini's theorem, we obtain,
\begin{multline} \label{Eq3DiscreteTnProof}
\int_{T_N\backslash
T_{N,j}}p(\mathbf{x})|u(\mathbf{x})|^2\mathrm{d}\mathbf{x}
 =\int_{\widehat{T}_N\backslash
\widehat{T}_{N,j}}\!\!p(\theta,\mathbf{s})|v(\theta,\mathbf{s})|^2J\mathrm{d}\xi\\[3mm]
\leq\!\!\!\int_{\widehat{T}_N\backslash
\widehat{T}_{N,j}}\!\!\frac{1}{C^2\hat{R}}\left(\int_{\theta}^{\hat{R}}\!
\left| \frac{\partial v}{\partial \vartheta}\right|^2
\,\mathrm{d}\vartheta\!\right)\!J\mathrm{d}\xi
 \leq \int_{\theta_j}^{\hat{R}}\!\!
\int_{\hat{\Omega}}\frac{1}{C^2\hat{R}}\left(\int_{\theta}^{\hat{R}}\!
\left| \frac{\partial v}{\partial \vartheta}\right|^2
\,\mathrm{d}\vartheta\!\right)\!\!J(\mathbf{s},\theta)\mathrm{d}\mathbf{s}\mathrm{d}\theta\\[3mm]
\leq\!\!\!\int_{\theta_j}^{\hat{R}}\!\!
\int_{\hat{\Omega}}\!\frac{1}{C\hat{R}}\!\left(\int_{\theta}^{\hat{R}}\!
\left| \frac{\partial v}{\partial
\vartheta}\right|^2\!\!J(\mathbf{s},\vartheta)
\,\mathrm{d}\vartheta\!\!\right)\!\!\mathrm{d}\mathbf{s}\mathrm{d}\theta
 \!\leq\!\frac{1}{C\hat{R}}\!\!\int_{\theta_j}^{\hat{R}}\!\!
\left(\int_{\hat{\Omega}}\int_{0}^{\hat{R}} \left| \frac{\partial
v}{\partial \vartheta}\right|^2\!\!
J(\mathbf{s},\vartheta)\mathrm{d}\vartheta\mathrm{d}\mathbf{s}\!\!\right)\!\!\mathrm{d}\theta
\\[3mm]
\leq \frac{1}{C}\int_{\widehat{T}_N}\!\! \left|\frac{\partial
v}{\partial \vartheta}\right|^2 J \mathrm{d}\xi
=\frac{1}{C}\int_{T_N}\!\! \left|\sum_{i=1}^n\frac{\partial
u}{\partial x_i}\frac{\partial x_i}{\partial \theta}\right|^2
\mathrm{d}\mathbf{x}
 \leq \int_{T_N} \!\!|\nabla u|^2
\mathrm{d}\mathbf{x}.
\end{multline}
Since \eqref{Eq1DiscreteTnProof} is satisfied, the spectrum of the
Laplacian on $T_N$ is purely discrete.  \qed
\begin{remark} {\em
A similar proof applies for a Schr\"{o}dinger operator on $T_N$
with a bounded from below potential.}\end{remark}

\subsection{Further results}
\label{SectionDiscApplications} In this subsection we present two
lemmas asserting that the $L^2$-norm of functions which are
bounded in $H^1_{0,\rho_\alpha}(T_1)$ and in $H_{0}^1(T_N^\eps)$
is concentrated on compact sets. These lemmas will be used in
Section~\ref{SectionPQ}. Their proofs are similar to those of
theorems~\ref{LemmaT1DiscLew} and \ref{lemmaDiscreteSpecT2}, and
therefore they  are omitted.

\begin{lemma} \label{LemmaL2NormNotInEndsT1}
    Assume that $T_1$ satisfies the assumptions of Theorem
    \ref{LemmaT1DiscLew}. Suppose that there exists a
    weight function $0<\rho<1$, which is constant on each edge of
    $T_1$, such that $\rho\asymp\rho_\alpha$ and $\rho\asymp\rho_\beta$
    with a constant $c$. Denote
    $T_{1,j}=\left\{\mathrm{gen}(e)\leq
    j\right\}$.
\begin{enumerate}
\item{
     Let $R(j)$ be the radius of the maximal (connected) subtree
    in $T\backslash T_{1,j}$. Then
    $$
        \int_{T_1\backslash T_{1,j}}  |u|^2 \rho_\beta\,\mathrm{d}\theta\leq
        \frac{c^2}{C}R(j)^2\int_{T_1}|u'|^2\rho_{\alpha}\,\mathrm{d}\theta
        \qquad \forall u\in H_{0,\rho_\alpha}^1(T_1).
    $$
    \label{LemmaT1Disc}
    }
\item{ Let $\overline{\mathcal{V}}^\eps:=\bigcup_{v\in T_1}
\overline{V}^\eps(v)$. Then
        $$
        \int_{\overline{\mathcal{V}}^\eps}  |u|^2 \rho_\beta\,\mathrm{d}\theta\leq
        O(\eps)\int_{T_1}|u'|^2\rho_{\alpha} \,\mathrm{d}\theta
       \qquad \forall u\in H_{0,\rho_\alpha}^1(T_1).
    $$ }
\end{enumerate}
\end{lemma}
\begin{lemma} \label{LemmaConvergenceAtEnds}
Assume that $T_N$ satisfies the assumptions
    of Theorem~\ref{lemmaDiscreteSpecT2}.
\begin{enumerate}
\item{ Let
    $T_{N,j}$ as defined in (\ref{Tnj}),
     and let $R(j):=\hat{R}-\theta_j$.
     Then
    \begin{equation}\label{DiscreteLemmaResult1}
        \int_{T_N\backslash T_{N,j}}|u(\mathbf{x})|^2\dxx\leq
        C^2R(j)^2\int_{T_N}|\nabla u|^2\dxx \qquad \forall u\in H_0^1(T_N).
    \end{equation}}
\item{ Let $\mathcal{V}^\eps:=\bigcup_{v\in T_1}
V^\eps(v)$. Then
    \begin{equation}
        \int_{\mathcal{V}^\eps}|u(\mathbf{x})|^2\dxx\leq
        O(\eps)\int_{T_N^\eps}|\nabla u|^2\dxx \qquad \forall u\in H_0^1(T_N^\eps).
    \end{equation}
 }
\end{enumerate}
\end{lemma}

\mysection{Convergence of the spectra of width-weighted operators}
\label{AnLimitSec} In this section we estimate the eigenvalues of
the width-weighted  operators on $T_1$ (defined in
Section~\ref{secDefA}), for the case where the weight functions
and the potential depend on $\eps$, and pointwise converge as
$\eps$ tends to $0$. We treat the weight functions and potential
term as convergent sequences of functions of $\eps$. Hence,
throughout this section we set $\eps:={1}/{n}$, where $n\in
\mathbb{N}$, and denote the weights and potentials by
$\rho_{\alpha,n},\rho_{\beta,n}$ and $W_n$. Accordingly, the
corresponding operators are denoted by $A_{\alpha,\beta,n}$, or
$A_n$ for short. We assume that $\rho_{\alpha,n}$ and
$\rho_{\beta,n}$ converge to a mutual weight function, which we
denote by $\rho$. We denote by $W$ the limit potential of the
sequence $W_n$. We also treat the spaces $\{H_{0,n}^1(T_1)\}:=
\{H_{0,\rho_{\beta,n}}^1(T_1)\}$ as a spaces sequence, with a
``limiting space" $H_{0,\rho}^1(T_1)$. Let
$\{L_n^2(T_1)\}_{n=1}^{\infty}$ and $L_\rho^2(T_1)$ be the
corresponding $L^2$ spaces. Using these notations, we study the
asymptotic behavior of the eigenvalues of $A_n$ as $n\rightarrow
\infty$.\\
Throughout this section we assume that the following conditions
are satisfied:
\begin{assumption} \label{AssumptionsWeightPotential}{\em
\begin{enumerate}
\item{$T_1$ has a finite radius.  }
\item{{\bf Assumptions on the weight functions:}
$\{\rho_{1,n}\}_{n=1}^{\infty}$ and
$\{\rho_{2,n}\}_{n=1}^{\infty}$ are positive bounded weight
functions sequences in $L^1_{\mathrm{loc}}(T_1)$, such that
$\rho_{1,n}\asymp\rho$ and $\rho_{2,n}\asymp\rho$ with the same
constant $c$ (so the spaces $H_{0,n}^1(T_1)$ and
$H_{0,\rho}^1(T_1)$ are equivalent for all $n\in \mathbb{N}$).
Moreover, for any neighborhood $U$ containing all the vertices of
$T_1$ and a given compact set $K\Subset T_1$, we have
$\rho_{1,n}=\rho_{2,n}= \rho$ in $(T_1\cap K)\setminus U$ for all
sufficiently large $n$.}
\item{{\bf Assumptions on the potential terms:}
$\{W_{T_1,n}\}_{n=1}^{\infty}$ is a sequence of real valued
radially symmetric potentials on $T_1$, for which there exists a
positive constant $C_W$ such that
$|W_{T_1,n}|_{L_\rho^\infty(T_1)}\leq C_W$. Moreover,
$\{W_{T_1,n}\}_{n=1}^{\infty}$ converges almost surely (and hence
in $L^1_{\rho,\mathrm{loc}}(T_1)$) to a potential $W$, which
satisfies $|W|_{L_\rho^\infty(T_1)}\leq C_W$. Without loss of
generality, we assume that $W_{T_1,n}>1$ for all $n\in
\mathbb{N}$. }
\end{enumerate}
 } \end{assumption}

Under Assumptions~\ref{AssumptionsWeightPotential}, we show that
the eigenvalues of the operators $A_n$ converge, as $n\rightarrow
\infty$,  to the eigenvalues of the limit operator $A$. Here the
operators $A_n$ are defined by the quadratic forms on
$H_{0}^1(T_1)\times H_{0}^1(T_1)$:
\begin{equation}\label{EqA_n}
    <A_n u,\phi>_n:=\int_{T_1}(u'\bar{\phi}'\rho_{1,n}
    +W_{T_1,n}u\bar{\phi}\,\rho_{2,n})\dt,
\end{equation}
while the limit operator $A$ is defined, similarly, by
\begin{equation}\label{EqA}
    <Au,\phi>:=\int_{T_1}(
    u'\bar{\phi}'+Wu {\bar{\phi}})\rho\dt.
   \end{equation}
 This result is stated in Corollary
\ref{MainTheorem}. Notice that since $\rho$ is constant on each
edge, the difference between the derivatives part of $A$ and the
Laplacian is manifested by the Kirchhoff condition.

In order to prove the convergence of the spectrum, we need the
following lemmas, whose proofs are given later.
\begin{lemma}\label{Lemma1}
    For $n\in \mathbb{N}$, consider operators  $A_n$  of the form \eqref{EqA_n} which
    satisfy Assumptions \ref{AssumptionsWeightPotential}.
    Let  $\{u_n\}\subset H_{0,\rho}^1(T_1)$ be a sequence of normalized eigenfunctions of $A_n$
    which converges weakly in $H_{0,\rho}^1(T_1)$
    to $u$.
        Let $\lambda_n$ be the sequence of corresponding eigenvalues of
    $A_n$. If
    $\lim_{n\rightarrow \infty} \lambda_{n}=\lambda$, then
$Au=\lambda u$, and $u\neq 0$ is an eigenfunction of the operator
$A$  defined in (\ref{EqA}) with eigenvalue $\lambda$. Moreover,
$\{u_n\}$ also converges locally uniformly to $u$.
\end{lemma}

\begin{lemma} \label{Lemma2}
    Consider operators $A$, and $A_n$ for $n\in \mathbb{N}$,   of the form \eqref{EqA_n} and
    \eqref{EqA} respectively which satisfy Assumptions
    \ref{AssumptionsWeightPotential}. Assume also that
    $A_n^{-1}:L^2_{\rho}(T_1)\rightarrow H^1_{0,\rho}(T_1)$ have
    uniform bounded norms.  Suppose that $u\neq 0$, and
    $Au=\lambda u$ in $L^2_\rho(T_1)$. For each $n\in \mathbb{N}$, let $w_n$  be
    the solution of the equation $A_nw_n=\lambda u$ in
    $L^2_{\rho}(T_1)$. Then $\{w_n\}$ has subsequence that we continue
    denoting by $\{w_n\}$, which converges to $u$ weakly in $H_{0,\rho}^1(T_1)$
    and strongly in $L^2_\rho(T_1)$. Moreover, $\{w_n\}$ also converges locally
    uniformly to $u$.
\end{lemma}

\begin{theorem}\label{MainTheorem}
    Let $\{\rho_{1,n}\}_{n=1}^{\infty}, \{\rho_{2,n}\}_{n=1}^{\infty}$
    and $\{W_{T_1,n}\}_{n=1}^{\infty}$ be sequences of weight functions and
    potentials on $T_1$ satisfying Assumptions
    \ref{AssumptionsWeightPotential}. Assume, in addition, that
    $\{W_{T_1,n}\}_{n=1}^{\infty}$ are continuous functions and that
     $\{\rho_{1,n}\}_{n=1}^{\infty}$ and $
    \{\rho_{2,n}\}_{n=1}^{\infty}$ equal $\rho$ except, at most, for
    $O(|e_j|/n)$ neighborhoods of vertices in generation $j$.
    Let a sequence of operators $A_n$ and a limit operator $A$ be
    defined by \eqref{EqA_n} and \eqref{EqA} respectively.
    We denote
    by $\lambda_{m,n}$ the $m$-th eigenvalue of $A_n$, and by
    $\lambda_m$ the $m$-th eigenvalue of $A$. Then
    $$\lim_{n\rightarrow \infty}\lambda_{m,n}=\lambda_m.$$
\end{theorem}
%
\begin{proof} We adapt Attouch's proof of \cite[Theorem~3.71]{Att}.
Since $\rho_{1,n}\asymp\rho$ and $\rho_{2,n}\asymp\rho$ with
a positive constant $c$, and $|W_{T_1,n}|$ and $|W|$ are bounded
by $C_W$, we have for all $u\neq 0$ that the Rayleigh quotients
satisfy
\begin{equation}\label{rq}R_n(u):= \frac{<\!A_nu,u\!>_{n}}{<\!u,u\!>_{n}}
=\frac{\int_{T_1}\!\!(|
u'|^2\rho_{1,n}+W_{T_1,n}|u|^2\rho_{2,n})\dtet}
{\int_{T_1}|u|^2\rho_{2,n}\dtet} \leq
c^2\frac{<\!Au,\!u>}{<\!u,u\!>_\rho}+2c^2C_W,
\end{equation}
and similarly
$$\frac{<A_nu,u>_{n}}{<u,u>_{n}}
\geq \frac{1}{c^2}\frac{<Au,u>}{<u,u>_\rho}-2\frac{1}{c^2}C_W.$$
Fix $l\in \mathbb{N}$, by the min-max principle we obtain
\begin{equation} \label{MainTheoremEq1}
    \frac{1}{c^2}(\lambda_l-2C_W)\leq \lambda_{l,n}\leq c^2(
    \lambda_l+2C_W),
\end{equation}
so, $\{\lambda_{l,n}\}$ is a bounded sequence. Therefore, there
exists a subsequence of $\{\lambda_{l,n}\}$ (that we keep denoting
by $\{\lambda_{l,n}\}$), and $\widehat{\lambda}_l\in \mathbb{R}$
such that  $\lambda_{l,n} \to \widehat{\lambda}_l$.

We claim that there exists an eigenfunction $\widehat{u}_l$ such
that $A \widehat{u}_l=\widehat{\lambda}_l \widehat{u}_l$,  i.e.,
$\{\widehat{\lambda}_l\}\subseteq\{\lambda_j\}$.
Indeed, let $\{u_{l,n}\}$ be the orthonormal sequence of
eigenfunctions of $A_n$ that correspond to $\{\lambda_{l,n}\}$. We
assume that $\|u_{l,n}\|_n=1$. Then
$$\int_{T_1}|\left( u_{l,n}\right)^{'}|^2\rho_{1,n}\,\mathrm{d}\theta
=\int_{T_1}(\lambda_{l,n}-W_{T_1,n})|u_{l,n}|^2\rho_{2,n}\,\mathrm{d}\theta
\leq\lambda_{l,n}+C_W.$$
 It follows that  $\{u_{l,n}\}$ is bounded in
$H^1_{0,\rho}$. The weak sequential compactness implies that
$\{u_{l,n}\}$ has a subsequence $\{u_{l,n}\}$ which converges
weakly in $H^1_{0,\rho}(T_1)$. We denote its limit by $
\widehat{u}_l$. By Lemma \ref{Lemma1}, $\widehat{u}_l\not= 0$,
$A\widehat{u}_l=\widehat{\lambda}_l\widehat{u}_l$ and the
convergence is locally uniform. In particular,
$\{\widehat{\lambda}_l\}\subseteq\{\lambda_j\}$. Moreover,
\eqref{MainTheoremEq1} implies that  $\{\widehat{\lambda}_l\}$ is
an infinite sequence, and since
$\{\widehat{\lambda}_l\}\subseteq\{\lambda_j\}$, we have  $\lim_{
l \rightarrow \infty} \widehat{\lambda}_l=\infty$.

Let us now show that $\{\lambda_j\} \subseteq
\{\widehat{\lambda}_l\}$. Assume that there exists an eigenvalue
$\lambda$ of $A$ such that $\lambda\neq \widehat{\lambda_l}$ for
all $l\in \mathbb{N}$, and let $u$ be a corresponding
eigenfunction of $A$ such that $||u||_{L_\rho^2}=1$.

Take  $m\in \mathbb{N}$ such that
$\lambda<\widehat{\lambda}_{m+1}$ for all limit values
$\widehat{\lambda}_{m+1}$ of the sequence $\{\lambda_{m+1,n}\}$.
Set $U_{m,n}= \mathrm{span}\{u_{1,n},...,u_{m,n}\}$. By the
min-max principle,
$$\lambda_{m+1,n}=\min_{v\in U_{m,n}^{\bot}}R_n(v),$$  where $R_n$ is defined in (\ref{rq}).
Therefore,  if we could find $v_n\in U_{m,n}^\bot$ satisfying
$$ \lim_{n\rightarrow\infty}R_n(v_n)\leq \lambda ,$$
then we would arrive to a contradiction of the assumption
$\lambda<\widehat{\lambda}_{m+1}$.

 Let  $w_n$ be the solutions
of the problem $A_nw_n=\lambda u$. The assumption $W_{T_1,n}>1$
implies that  $A_n^{-1}$ are uniformly bounded, so $\{w_n\}$ is a
bounded sequence in $L^2_{\rho}(T_1)$. By Lemma \ref{Lemma2}, up
to a subsequence, $\{w_n\}$ converges to $u$, weakly in
$H_{0,\rho}^1(T_1)$, strongly in $L^2_\rho(T_1)$, and also locally
uniformly.

Let us show that $\lim_{n\rightarrow \infty}R_n(w_n)=\lambda:$
\begin{multline} \label{Eq1ProofTheoremConve}
<A_nw_n,w_n>_n=\lambda<u,w_n>_n
=\lambda\int_{T_1}\!\![u\overline{w}_n(\rho_{2,n}-\rho)+u(\overline{w}_n-\overline{u})\rho]\dt
+\lambda.
\end{multline}
Since
$$ \left|\int_{T_1}\!\!\!u\overline{w}_n(\rho_{2,n}-\rho)\dt\right|^2
\!\!\!=\left|\int_{T_1}\!\!\!u\overline{w}_n\rho\frac{(\rho_{2,n}-\rho)}{\rho}\dt\right|^2
\leq \left\|
u\left(\frac{\rho_{2,n}-\rho}{\rho}\right)\right\|^2_{L^2_\rho}
\left\|
 w_n\right\|^2_{L^2_\rho} \ ,
$$
Lebesgue's dominated convergence theorem implies that the first
term of the right-hand side of \eqref{Eq1ProofTheoremConve}
converges to zero, while the second term tends to zero due to the
$L^2_\rho(T_1)$ convergence of $\{w_n\}$ to $u$.
Therefore,
\begin{equation}\label{Eq2ProofTheoremConve}
    \lim_{n\rightarrow \infty}<A_nw_n,w_n>_n=\lambda.
\end{equation}
Moreover,
\begin{equation}\label{Eq3ProofTheoremConve}<w_n,w_n>_n
=\int_{T_1}\!\!\left[(|w_n|^2-|u|^2)\rho_{2,n}
+|u|^2(\rho_{2,n}-\rho)+|u|^2\rho\right]\dt.
\end{equation}
The first terms in (\ref{Eq3ProofTheoremConve}) converges to zero
due to  the strong convergence of $w_n$ to $u$ in
$L^2_{\rho}(T_1)$. Indeed,
$$ \left|\int_{T_1}\!\!\!\!\left(|w_n|^2\!-\!|u|^2\right)\rho_{2,n}\dt\right|
  \leq  \left\| w_n\!-\!u\right\|^{1/2}_{
{L}^2_{\rho_n}} \!\left\| w_n\!+\!u\right\|^{1/2}_{
L^2_{\rho_n}}\leq C
 \left\| w_n\!-\!u\right\|^{1/2}_{
{L}^2_{\rho}}\! \left\| w_n\!+\!u\right\|^{1/2}_{ L^2_{\rho}}.  $$
The second term in (\ref{Eq3ProofTheoremConve}) converges to zero
by Lebesgue's dominated convergence theorem.
Hence, \eqref{Eq2ProofTheoremConve} and
\eqref{Eq3ProofTheoremConve} imply that
\begin{equation} \label{Eq5ProofTheoremConve}
    \lim_{n\rightarrow\infty}R_n(w_n)=\lambda.
\end{equation}
Define
$$v_n:=w_n-\sum_{k=1}^m<w_n,u_{k,n}>_n u_{k,n}.$$
Fix $1\leq k\leq m$, and let $\widehat{u}_k$ be a weak limit of
$u_{k,n}$. It follows (as above) that
\begin{equation} \label{Eq5_5ProofTheoremConve}
    \lim_{n\rightarrow\infty}<w_n,u_{k,n}>_n =
    \lim_{n\rightarrow\infty}<w_n,u_{k,n}>_\rho
    =<u,\widehat{u}_k>_\rho.
\end{equation}
 By the first part of the proof,
$\widehat{u}_k$ is an eigenfunction of $A$, and by our assumption,
its eigenvalue is not equal to $\lambda$. Therefore,
$<u,\widehat{u}_k>_n=0$ and by \eqref{Eq5_5ProofTheoremConve},
\begin{equation} \label{Eq6ProofTheoremConve}
    \lim_{n\rightarrow \infty} <w_n,u_{k,n}>_n=0.
\end{equation}
That implies that $\{v_n\}$ and $\{w_n\}$ share the same
$L^2$-limit $u$.

Using \eqref{Eq2ProofTheoremConve} and
\eqref{Eq6ProofTheoremConve}, a direct calculation yields that
$$
\lim_{n\rightarrow \infty}\!\!<\!A_nv_n,v_n\!>_n =\lambda,
\>\>\mbox{  and  }\>\> \lim_{n\rightarrow \infty}<v_n,v_n>_n =1.
$$
Hence
$$\lim_{n\rightarrow \infty} R_n(v_n)
=\lim_{n\rightarrow \infty} R_n(w_n) =\lambda.$$ By the definition
of $v_n$, we have $<v_n,u_{k,n}>_{L^n_n}=0$ for all $k=1,...,m$.
Hence, the min-max principle implies that $R_n(v_n)\geq
\lambda_{m+1,n}$. Therefore, $\lambda\geq \widehat{\lambda}_{m+1}$
for some limit value $\widehat{\lambda}_{m+1}$, which contradicts
the assumption $\lambda<\min \{\widehat{\lambda}_{m+1}\}$.
\end{proof}

\begin{remark}\label{LemmaArzela} {\rm
%
Let $\{u_n\}_{n=1}^{\infty}\subseteq H^1_{0,\rho}(T_1)\cap
C^1(T_1)$ be a sequence which converges weakly to $u$ in
$H^1_{0,\rho}(T_1)$. It follows that $\{u_n\}$ is locally a
bounded and equicontinuous sequence in $C(T_1)$. By
Arzel\`{a}-Ascoli's theorem, $\{u_n\}_{n=1}^{\infty}$ has a
subsequence that converges locally uniformly to a continuous
function $u$.
}
\end{remark}

\begin{proof}[{Proof of Lemma \ref{Lemma1}}]  By
Remark~\ref{LemmaArzela}, $\{u_n\}$ has a subsequence which we
continue denoting by $\{u_n\}$, that converges locally uniformly
to $u$ which is continuous on $T_1$. We claim:  (1) $u\in
\mathrm{Dom}(A)$,  (2) $Au=\lambda u$ and , and (3) $u\neq 0$. The
first two claims follow provided we prove
\begin{equation}\label{Eqr0.5}
    \int_{T_1}(u'\overline{\phi'}+W u\overline{\phi})\rho\dtet
    =\lambda \int_{T_1}  u\overline{\phi}\rho \dtet \qquad \forall \phi\in C_0^1(T_1).
\end{equation}
Since $\{u_n\}$ are eigenfunctions of $A_n$, for each test
function $\phi\in C_0^1(T_1)$,
\begin{equation} \label{Eqr1}
    \int_{T_1}(u_n'\overline{\phi'}\rho_{1,n}+W_{T_1,n}u_n\overline{\phi}\rho_{2,n})\dtet
    =\int_{T_1}\lambda_n u_n \overline{\phi} \rho_{2,n}\dtet.
\end{equation}
By Lebesgue's theorem applied to $\rho_{1,n}$ and the
$H^1_{0,\rho}$ bound of $u_n$,
\begin{equation} \label{Eqr2}
    \lim_{n\rightarrow\infty}\left|
    \int_{T_1}\!\!u_n' \overline{\phi'}(\rho-\rho_{1,n})\dtet\right|^2 \leq
    \lim_{n\rightarrow\infty}
    \int_{T_1}\!\!\!|\phi'|^2\frac{(\rho-\rho_{1,n})^2}{\rho}\dtet
    \int_{T_1} |u_n'|^2\rho\dtet\!=\!0.
\end{equation}

The weak convergence of $\{u_n\}$ to $u$ in $H^1_{0,\rho}(T_1)$
implies that
\begin{equation}\label{Eqr2.5}
    \lim_{n\rightarrow\infty }\int_{T_1}
    u_n'\overline{\phi'}\rho\dtet
    =\int_{T_1} u'\overline{\phi'}\rho\dtet.
\end{equation}

By similar  arguments, the local uniform convergence of $\{u_n\}$
to $u$, and the a.s. convergence of $W_{T_1,n}$ and $\rho_{2,n}$
imply that
\begin{equation}\label{Eqr2.6} \lim_{n\rightarrow \infty}\lambda_n\!\!\int_{T_1}\!\!
u_n \overline{\phi}\rho_{2,n}\dtet\!=\!\lambda\!\!\int_{T_1}\!\! u
\overline{\phi}\rho\dtet, \; \mbox{and }
\lim_{n\rightarrow\infty}\!\!\int_{T_1}\!\!
W_{T_1,n}u_n\overline{\phi}\rho_{2,n}\dtet\!=\! \int_{T_1}\!\!
Wu\overline{\phi}\rho\dtet.
\end{equation}
Now, \eqref{Eqr1}--\eqref{Eqr2.6} imply \eqref{Eqr0.5}.

In order to show that $u\neq 0$, let $T_{1,k}=\{e\in
T_1\>|\>\mathrm{gen}(e)\leq k\}$, and let $R(k)$ be the maximal
radius of subtrees in $T_1\backslash T_{1,k}$.  Recall that
$\{u_n\}$ are eigenfunctions satisfying
$\|u_n\|_{H^1_{0,\rho}(T_1)}=1$, the corresponding eigenvalues
sequence $\{\lambda_n\}$ converges, the potential terms $\{W_n\}$
are bounded by a constant $C_W$ for all $n\in \mathbb{N}$, and
$\rho_{1,n}\asymp\rho\asymp \rho_{2,n}$. Therefore, using
\eqref{Eqr2} and the arguments that eigenfunctions has
$L^2_\rho(T_1)$ and $H^1_\rho(T_1)$ norms of the same order, we
infer that there exist $\gamma, \delta>0$ so that,  for $n$ large
enough,  $||u_n||_{L^2_\rho(T_1)}\geq \gamma>0\>\>\mbox{ and }\>\>
||u_n'||_{L^2_\rho(T_1)}\leq \delta$.

Therefore, by Lemma~\ref{LemmaL2NormNotInEndsT1}
(\ref{LemmaT1Disc}) we have that
\begin{align*}
\int_{T_{1,k}} |u_n|^2\rho\,\mathrm{d}\theta =\int_{T_1}
|u_n|^2\rho\,\mathrm{d}\theta-\int_{T_1\backslash T_{1,k}}
|u_n|^2\rho\,\mathrm{d}\theta \\[2mm]
 \geq \int_{T_1}
|u_n|^2\rho\,\mathrm{d}\theta-\frac{c^2}{C}R(k)^2\int_{T_1}
|u_n'|^2\rho\,\mathrm{d}\theta \geq
\gamma-\delta\frac{c^2}{C}R(k)^2.
\end{align*}
Now, choose $k$ large enough such that
$\gamma-\delta[cR(k)]^2/C>0$. By the local uniform convergence of
$u_n$ to $u$, we obtain
$$0<\gamma-\delta\frac{c^2}{C}R(k)
\leq \lim_{n\rightarrow \infty}\int_{T_{1,k}}
|u_n|^2\rho\mathrm{d}\theta
=\int_{T_{1,k}} |u|^2\rho\mathrm{d}\theta
\leq \int_{T_1} |u|^2\rho\mathrm{d}\theta.
$$ Therefore $u\neq 0$, and $u$ is an eigenfunction of $A$.
\end{proof}


\begin{proof}[{Proof of Lemma \ref{Lemma2}}]
Since $Au=\lambda u$ and $A$ is invertible, it is sufficient to
prove that $Aw=Av$ (and in particular that $w$ is in the domain of
$A$). But this is equivalent to \begin{equation}\label{ww3}
\left\langle A w, \phi\right\rangle = \left\langle A u,
\phi\right\rangle
\end{equation}
for any function $\phi$ in a dense subset of $H^1_{0,\rho}(T_1)$.
Recall that $w\in H^1_{0,\rho}(T_1)$ and $\left\langle A w,
\phi\right\rangle$ is defined by (\ref{EqA}). Let us split the
quadratic form (\ref{EqA}) into
$$ \left\langle A w,
\phi\right\rangle  = \left\langle A w, \phi\right\rangle^{(1)} +
\left\langle A w, \phi\right\rangle^{(2)},$$ where
\begin{equation}\label{ww1} \left\langle A w, \phi\right\rangle^{(1)}:=
\int_{T_1}
    w'\bar{\phi}'\rho \dt, \ \ \  \left\langle A w,
    \phi\right\rangle^{(2)}:=\int_{T_1}Ww
    {\bar{\phi}}\rho\dt .
\end{equation} Similarly, (\ref{EqA_n}) is written as
$$ \left\langle A_n w,
\phi\right\rangle_n  = \left\langle A_n w,
\phi\right\rangle_n^{(1)} + \left\langle A_n w,
\phi\right\rangle_n^{(2)},$$ where \begin{equation}\label{ww2}
\left\langle A_n w, \phi\right\rangle_n^{(1)}:= \int_{T_1}
    w'\bar{\phi}'\rho_{1,n} \dt, \ \ \  \left\langle A w,
    \phi\right\rangle^{(2)}_n:=\int_{T_1}Ww
    {\bar{\phi}}\rho_{2,n}\dt .
\end{equation}

Let $\Phi(T_1)$ be  the set of all functions $\phi\in C^2_0(T_1)$
which are constant in some neighborhood of any vertex $v\in T_1$.

 We further
observe
\begin{enumerate}
  \item  For any  $\phi\in\Phi(T_1)$ and sufficiently large $n$, $\phi^{'}=0$
whenever  $\rho_{1,n}\not=\rho$ or
  $\rho_{2,n}\not=\rho$ by Assumption~\ref{AssumptionsWeightPotential} (2).
Hence, for  a given  $\phi\in\Phi(T_1)$
$$ \langle A_nw_n, \phi\rangle^{(1)}_{n}=\langle
A w_n,\phi\rangle^{(1)} .
$$
for all sufficiently large $n$.
\item Since $w$ is the weak
limit of $w_n$ in $H^1_{0,\rho}(T_1)$ it follows by (1)  that
$$ \lim_{n\rightarrow\infty}\langle
A_n w_n,\phi\rangle^{(1)}_n =\lim_{n\rightarrow\infty}\langle A
w_n,\phi\rangle^{(1)}=\langle A w ,\phi\rangle^{(1)} .
$$
\item By Assumption~\ref{AssumptionsWeightPotential} and the
strong convergence of $w_n$ to $w$ in $ L^2_\rho(T_1)$ we obtain
$$ \lim_{n\rightarrow\infty}
\left\langle A_n w_n, \phi\right\rangle^{(2)}_n= \left\langle A w,
\phi\right\rangle^{(2)}  . $$
 \item By (2) and (3) we obtain
$$ \lim_{n\rightarrow\infty}
\left\langle A_n w_n, \phi\right\rangle_n= \left\langle A w,
\phi\right\rangle  . $$ for any $\phi\in \Phi(T_1)$.
 \item
Since $A_nw_n=\lambda u=Au$ by assumption we obtain $\left\langle
A_n w_n, \phi\right\rangle_n= \lambda\left\langle u,
\phi\right\rangle_n= \left\langle A u, \phi\right\rangle_n $.
Since $\rho_{n,1}\rightarrow \rho$ in measure, it follows that
$$ \left\langle A w,
\phi\right\rangle=\lim_{n\rightarrow\infty} \left\langle A_n w_n,
\phi\right\rangle_n= \lambda\left\langle u, \phi\right\rangle_\rho
= \left\langle A u, \phi\right\rangle  .
$$
\end{enumerate}
So,  (\ref{ww3}) is proved for any $\phi\in\Phi(T_1)$. The proof
is completed by observing that $\Phi(T_1)$ is clearly dense in
$H^1_{0,\rho}(T_1)$.
\end{proof}

\mysection{$\eps$-dependent bounds for the eigenvalues of
$N$-dimensional tree} \label{SectionPQ}

In this section, we consider the spectrum of the Schr\"{o}dinger
operator $$L_\eps :=-\Delta +W_{T_N^\eps}$$ on $T_N^\eps$, where
$N$ is dimension of the tree, and $W_{T_N^\eps}$ is a continuous
bounded potential on $T_N^\eps$. Without loss of generality, we
assume that $W_{T_N^\eps}\geq 0$.

We prove that the eigenvalues of $L_\eps$ are bounded from above
and below by functions $\phi_Q^\eps,\phi_P^\eps$ of the
eigenvalues of weighted operators $A_\eps$ on $T_1$ of the form
\begin{equation} \label{DefWidthWeightedA}
    A_\eps:=-\frac{1}{\rho_{b,\eps}}
    \frac{\mathrm{d}}{\mathrm{d}\theta}\left(\rho_{a,\eps}
    \frac{\mathrm{d}}{\mathrm{d}\theta}\right)+W_{T_1,\eps}\,,
\end{equation}
for a suitable choice of weight functions $\rho_{a,\eps},
\rho_{b,\eps}$ and a potential $W_{T_1,\eps}$ on $T_1$ of the form
\begin{equation} \label{DefV_t1}
W_{T_1,\eps}(\theta):=\left\{%
\begin{array}{ll}
    \dfrac{\int_{\Omega^\eps_e}W_{T_N^\eps}(\theta,\mathbf{s})\,\mathrm{d}\mathbf{s}}
    {|\Omega^\eps_e|}
            & \theta\in \overline{E}^\eps(e), \\[5mm]
    \sum_{e\in N(v)}b_e\psi_{(e)}(\theta) & \theta\in \overline{V}^\eps(v),  \\
\end{array}%
\right.
\end{equation}
where $b_e=(|\Omega^\eps_e|)^{-1}\int_{\Omega_e^\eps}
W_{T_N^\eps}(p_e^\eps,\mathbf{s})\,\mathrm{d}\mathbf{s}$,
$p_e^\eps=\partial\overline{V}^\eps(v)\cap e$ is the end point of
$\overline{V}^\eps(v)$ corresponding to $e\in N(v)$, and
$\{\psi_{(e)}\}$ is the partition of unity in a neighborhood of
the vertex $v$ defined in Section~\ref{appen1}. The functions
$\phi_Q^\eps$ and $\phi_P^\eps$ converge to the identity function
as $\eps$ tends to zero.
\subsection{Rayleigh quotients of Schr\"{o}dinger operator on $T_1$
and $T_N^\eps$}

The comparison between the Rayleigh quotients on $T_1$ and
$T_N^\eps$ involves the construction of transformations
$Q^\eps:H^1_{0,\rho^*}(T_1)\rightarrow H_0^1(T_N^\eps)$ and
$P^\eps:H_0^1(T_N^\eps)\rightarrow H^1_{0,\rho^*}(T_1)$, where
$\rho^*:T_1\rightarrow \mathbb{R}$ is defined by
$\rho^*(\theta):=\delta_e^{(N-1)}|\Omega_e|$ for $\theta\in e$. We
devote the following two subsections for the definitions of these
transformations.

\subsubsection{The mapping $Q^\eps:H^1_{0,\rho^*}(T_1)
\rightarrow H_0^1(T_N^\eps)$.} \label{QSection} Given a function
$f\in H^1_{\rho^*}(T_1)$ and a vertex $v$, we denote by
$f^\eps_e=f(p_e^\eps)$
 and
$\vec{f}^\eps=\left\{f_e^\eps\right\}_{e\in N(v)}$. $Q^\eps(f)$ is
defined as follows:
\begin{equation} \label{QDef}
Q^\eps(f)(\mathbf{x})=
\left\{%
\begin{array}{ll}
    f(\theta) & \mathbf{x}=(\theta,s)\in E^\eps,\\
    \sum_{e\in N(v)} f_e^\eps\phi_{(e)}^\eps & \mathbf{x}\in V^\eps(v) , \\
\end{array}%
\right.
\end{equation}
where  $\{\phi^\eps_{(e)}\}$ is a partition of unity of
$V^\eps(v)$ as defined in Section~\ref{appen1}. We denote
$Q(f):=Q^1(f)$. We also define
\begin{equation} \label{DefrhoQ}
\rho_{Q}^\eps(\theta)=\left\{%
\begin{array}{ll}
    \rho^* & \theta\in \overline{E}^\eps(e), \\
    \max\left\{\frac{\alpha^A}{\beta^{\overline{A}}},
    \frac{\alpha^B}{\beta^{\overline{B}}}\right\}\rho^*
        & \theta\in \overline{V}^\eps(v), \\
\end{array}%
\right.
\end{equation}
where $\alpha^A,\beta^{\overline{A}},\alpha^B$ and
$\beta^{\overline{B}}$ are defined in Section~\ref{appen1} and
Lemma~\ref{LemmaAB}.

\begin{lemma}\label{LemmaQ} There exists $c>0$ such that for any  $f\in H^1_{0,\rho^*}(T_1)$
and $0<\eps<1$, we have


\begin{enumerate}
\item{$Q^\eps(f)\in L^2(T_N^\eps)$.}
\vspace {2mm}

\item{$\int_{T_N^{\eps}}|\nabla Q^\eps(f)|^2\dxx
\leq \eps^{(N-1)}\!\int_{T_1}|f^{'}|^2\rho^\eps_Q\dtet \ $.
Moreover, $Q^\eps(f)\in H_0^1(T_N^\eps)$. }

\vspace {2mm}

\item{
$\int_{T_N^\eps}\!|Q^\eps(f)|^2\!\dxx\!\geq\!
\eps^{(N-1)}\!\int_{T_1}\!\!(|f|^2\!-\!c\eps|f^{'}|^2)\rho^*\!\dtet
$.}
\vspace {2mm}

\item{$\int_{T_N^\eps}W_{T_{n}^{\eps}}
|Q^\eps(f)|^2\dxx \leq \eps^{(N-1)}\!\int_{T_1}
(W_{T_1,\eps}|f|^2+O(\eps)|f^{'}|^2)\rho^*\dtet$. }
\end{enumerate}
\end{lemma}

\begin{proof} 1. We denote a
normalized connector $\overline{V}:=
(\eps\delta)^{-1}\overline{V}^\eps(v)$, where
$\delta=\delta_v=\delta^{\mathrm{gen}(v)}$ corresponds to the
vertex $v$ in question, and
\begin{equation}\label{inq}\widehat{f}(\theta):=f\left(\frac{\theta}{\eps
\delta}\right),\qquad
\widehat{\rho^*}(\theta):=\rho^*\left(\frac{\theta}{\eps
\delta}\right)
\end{equation} are the representation of $f$ and $\rho^*$ in $\overline{V}$
(here $\theta=0$ corresponds to $v$).
$$\int_{E^\eps}\!\!\! |Q^\eps(f)|^2\dxx\!=\!
\int_{\overline{E}^\eps}\! \int_{\Omega^\eps}\!\! |f(\theta)|^2
\mathrm{d}\mathbf{s}\mathrm{d}\theta
\!=\!\int_{\overline{E}^\eps}\!\!(\eps\delta)^{(N-1)}
|\Omega^\eps||f|^2\dtet
\!=\! \eps^{(N-1)}\!\! \int_{\overline{E}^\eps}\!|f|^2\rho^* \dtet
.$$
For the connector $V^\eps$, we have by Lemma~\ref{LemmaAB}
 that
$$ \int_{V^\eps} |Q^\eps(f)|^2\dxx
= \!=\!(\eps \delta)^N\!\!\! \int_{V}\!\!\!|Q(f)|^2 \dxx\leq (\eps
\delta)^N\alpha^B|{\vec{f}}|^2.
$$
By Lemma~\ref{LemmaFormOnVBarj} and  since $\delta<1$, we have
that
$$
(\eps \delta)^N\alpha^B|{\vec{f}}|^2
\leq \frac{(\eps
\delta)^N\alpha^B}{\delta^{(N-1)}\beta^{\overline{B}}}
\int_{\overline{V}}\!\!\!(|\widehat{f}'|^2+|\widehat{f}|^2)\widehat{\rho}^*\dtet
 \leq\frac{\eps^{(N-1)}\alpha^B}{\beta^{\overline{B}}}
\int_{\overline{V}^\eps}\!\!\!(|\eps f'|^2+ |f|^2)\rho^* \dtet .
$$
In particular, we proved
\begin{equation}\label{inpar}\|Q^\eps
(f)\|_{L^2(T_N^\eps)}^2 \leq
\eps^{(N+1)}\|f'\|^2_{L^2_{\rho_Q^\eps}(T_1)}+
\eps^{(N-1)}\|f\|^2_{L^2_{\rho_Q^\eps}(T_1)}  .
\end{equation}

2. \
$$\int_{E^\eps} |\nabla
Q^\eps(f)|^2\dxx =\int_{E^\eps} |f^{'}|^2\dxx
=\int_{\overline{E}^\eps}\int_{ \Omega^\eps} |f^{'}|^2
\mathrm{d}\mathbf{s}\mathrm{d}\theta =\eps
^{(N-1)}\int_{\overline{E}^\eps} |f^{'}|^2 \rho^*\dtet .
$$

For  $V^\eps$ we use similar considerations to those used in part
1:
\begin{multline*}
\int_{V^\eps}\!\!\! |\nabla Q^\eps (f)|^2\dxx
=(\eps\delta)^{(N-2)}\int_{V} \!\!|\nabla Q (f)|^2\dxx
=(\eps\delta)^{(N-2)}{\vec{f}} A{\vec{f}}^* \\\leq\!
(\eps\delta)^{(N-2)}\alpha^A|{\vec{f}}\llcorner {\vec{1}}|^2
\!\leq\!\frac{(\eps\delta)^{(N-2)}}{\delta^{(N-1)}}\frac{\alpha^A}{\beta^{\overline{A}}}
\!\int_{\overline{V}}\!|\widehat{f}'|^2\rho^*\!\dtet
\!=\!\eps^{(N-1)}\frac{\alpha^A}{\beta^{\overline{A}}}
\!\int_{\overline{V}^\eps}\!|f'|^2\rho^*\!\dtet .
\end{multline*}
So, $Q^\eps(f)\in H_0^1(T_N^\eps)$ by definition (\ref{DefrhoQ}).

3. By Lemma~\ref{LemmaL2NormNotInEndsT1},
\begin{multline*}
\int_{T_N^\eps} |Q^\eps(f)|^2\dxx \geq
\int_{T_N^\eps\setminus\cup_{v}V^\eps(v)} |Q^\eps(f)|^2\dxx
=\eps^{(N-1)}\int_{\cup_{e}\overline{E}(e)}|f|^2\rho^*\dtet \\
\!=\!\eps^{(N-1)}\!\!\int_{T_1}\!|f|^2\rho^*\dtet\!-\!
\eps^{(N-1)}\!\!\int_{\cup_{v}\overline{V}^\eps(v)}\!\!|f|^2\rho^*\!\dtet\!
\geq\! \eps^{(N-1)}\!\!\int_{T_1}\!\!|f|^2\rho^*\!\dtet
\!-\!c\eps^{N}\!\!\int_{T_1}\!\!|f'|^2\rho^*\!\dtet .
\end{multline*}
The proof of (4) is  a simple extension of (\ref{inpar}).
\end{proof}
\begin{corollary} \label{Q_eps_Theorem}
There exists a constant $c>0$ such that for all $f\in
H^1_{0,\rho^*}(T_1)$ and $0<\varepsilon$ sufficiently small, the
Rayleigh quotients
$$R_{H_0^1(T_N^\eps)}[Q^\eps f]
:=\frac{\int_{T_N^\eps}(|\nabla (Q^\eps
f)|^2\!+\!W_{T_{n}^{\eps}}|Q^\eps(f)|^2)
\dxx}{\int_{T_N^\eps}|Q^\eps f|^2\dxx}\,,
$$
and
$$
R_{H^1_{0,\rho^*}(T_1)}[f]
:=\frac{\int_{T_1}(|f_\theta|^2\rho_Q^\eps\!+\!W_{T_1,\eps}|f|^2\rho^*)\dtet
}{\int_{T_1} |f|^2\rho^*\dtet}
$$ satisfy the inequality
\begin{equation}\label{ReqSat1_0}
     R_{H_0^1(T_N^\eps)}[Q^\eps f]
     \leq\frac{(1+O(\eps))R_{H^1_{0,\rho^*}(T_1)}[f]}
     {1-c\eps R_{H^1_{0,\rho^*}(T_1)}[f]}\,.
\end{equation}
\end{corollary}

\begin{remark} {\em Notice that
$R_{H^1_{0,\rho^*}(T_1)}[f]$ depends on $\eps$ and is the Rayleigh
quotient of the width-weighted operator $A_\eps$ defined in
\eqref{DefWidthWeightedA}, substituting
$\rho_{\alpha,\eps}=\rho_Q^\eps$ and $\rho_{b,\eps}=\rho^*$. }
\end{remark}

\subsubsection{The mapping $P^\eps:
H_0^1(T_N^\eps)\rightarrow H^1_{0,\rho^*}(T_1)$} \label{PSection}

Given a function $u\in H^1_0(T_N^\eps)$, a vertex $v$ and edges
$e\in N(v)$,  we denote
$$u_e:=\frac{1}{|\Omega_e^\eps|}
\int_{\Omega_e^\eps}u(p_e^\eps,\mathbf{s})\,\mathrm{d}\mathbf{s},
 \>\>\> \>\>\> \>\>\>
\vec{u}={\vec{u}_v}:=\left\{u_e\right\}_{e\in V(v)},$$  where
$p_e^\eps=\partial\overline{V}^\eps(v)\cap e$ are the end points
of $\overline{V}^\eps(v)$. Define
\begin{equation} \label{Pdef}
P^\eps(u)(\theta):=
\left\{%
\begin{array}{ll}
    \dfrac{\int_{\Omega^\eps}u(\theta,\mathbf{s})\mathrm{d}\mathbf{s}}{|\Omega^\eps|}
        & \theta\in \overline{E}^\eps,\\[4mm]
    \displaystyle\sum_{e\in V(v)} u_e\psi_{(e)}^\eps(\theta) & \theta\in \overline{V}^\eps , \\
\end{array}%
\right.
\end{equation}
where $\{\psi_{(e)}\}$ is the partition of unity in a neighborhood
of the vertex $v$ defined in Section~\ref{appen1}.
We also define
\begin{equation}\label{DefrhoP}\rho_{P}^\eps(\theta):=\left\{%
\begin{array}{ll}
    \rho^* & \theta\in \overline{E}^\eps , \\
    \min\left\{\frac{\beta^A}{\alpha^{\overline{A}}},
    \frac{\beta^B}{\alpha^{\overline{B}}}\right\}\rho^*
       & \theta\in \overline{V}^\eps. \\
\end{array}%
\right.
\end{equation}

\begin{lemma} \label{LemmaP}
There exists $c>$ such that for any  $u\in H^1_0(T_N^\eps)$ and
$0<\eps$ sufficiently small, we have


\begin{enumerate}
\item $P^\eps(u)\in L^2_{\rho^*}(T_1)$.

\vspace {4mm}

\item $\eps^{(N-1)}\int_{T_1}|(
P^\eps u)^{'}|^2\rho_P^\eps\dtet \leq \int_{T_N^\eps}|\nabla u|^2
\dxx$. In particular, $P^\eps(u)\in H_{0,\rho^*}^1(T_1)$.

\vspace{4mm}

\item $\eps^{(N-1)}\int_{T_1}|P^\eps
u|^2\rho^*\dtet \geq \int_{T_N^\eps}[(1-\sqrt{\eps})|u|^2-c\eps
|\nabla u|^2]\dxx$.

\vspace{4mm}

\item  $\eps^{(N-1)}\int_{T_1} W_{T_1,\eps}|P^\eps
u|^2\rho^*\dtet \leq
\int_{T_N^\eps}(1+2\sqrt{\eps})(W_{T_N^\eps}|u|^2+O(\eps) |\nabla
u|^2) \dxx $.
\end{enumerate}
\end{lemma}
\begin{proof}
Throughout the proof we denote
$\widehat{u}(\mathbf{x})=u\left(\delta\eps\mathbf{x}\right)$ for
$\mathbf{x}\in V$ (so, $\eps\delta\mathbf{x}\in V^\eps$).
Similarly,
$\widetilde{u}(\theta,\mathbf{s})=u\left(\theta,\eps\delta
\mathbf{s}\right)$ for $( \theta,\mathbf{s})\in
\overline{E}\times\Omega$ (so, $(\theta,\eps\delta\mathbf{s})\in
E^\eps$).

1. 
For each edge $E^\eps$,
\begin{multline*}
\eps^{(N-1)} \int_{\overline{E}^\eps}|P^\eps u|^2\rho^*\dtet
=\int_{\overline{E}^\eps}|\Omega^\eps||P^\eps u|^2\dtet
=\int_{\overline{E}^\eps}|\Omega^\eps|
\left|\frac{1}{|\Omega^\eps|}
\int_{\Omega^\eps}u(\theta,\mathbf{s})\mathrm{d}\mathbf{s}\right|^2\mathrm{d}\theta
\\[2mm]
\leq \int_{\overline{E}^\eps}
\int_{\Omega^\eps}|u(\theta,\mathbf{s})|^2\mathrm{d}\mathbf{s}\mathrm{d}\theta
=\int_{E^\eps}|u(\theta,\mathbf{s})|^2\mathrm{d}\mathbf{x}.
\end{multline*}
Using Lemma~\ref{LemmaAB}, we obtain for the connector $V^\eps$
$$\eps^{(N-1)}\!\!\!\int_{\overline{V}^\eps}\!\!|P^\eps
 u|^2\rho^*\!\dtet\!
=\!\eps^{N}\delta\!\!\sum_{e,\tilde{e}\in
N(v)}u_e\overline{u_{\tilde{e}}}
\!\!\int_{\overline{V}}\!\psi_{(e)}\psi_{(\tilde{e})}\rho^*\!\dtet
\!=\!\eps^{N}\delta {\vec{u}}\overline{B}(\vec{u})^*
\!\leq\!(\eps\delta)^N \alpha^{\overline{B}}|{\vec{u}}|^2.
$$
By Lemma~\ref{LemmaFormOnVj}, and assuming $\delta<1$, we obtain
$$
(\eps\delta)^N \alpha^{\overline{B}}|{\vec{u}}|^2 \leq
(\eps\delta)^{N}
\frac{\alpha^{\overline{B}}}{\beta^B}\int_{V}(|\nabla
\widehat{u}|^2+|\widehat{u}|^2)\dxx
\leq\frac{\alpha^{\overline{B}}}{\beta^B}\int_{V^\eps}
(|\eps\nabla u|^2+|u|^2)\dxx .
$$

2. For an edge $E^\eps$ we have
\begin{multline}\label{eqProofPeps0_5}
    \eps^{(N-1)}\!\!\int_{\overline{E}^\eps}
    \left|\left( P^\eps u\right)^{'}\right|^2\!\!\rho^*\dtet
    =\eps^{(N-1)}\!\!\int_{\overline{E}^\eps}\left|
    \frac{1}{|\Omega^\eps|}\int_{\Omega^\eps}\frac{\partial
    u(\theta,s)}{\partial
    \theta}\mathrm{d}\mathbf{s}\right|^2\!\!\rho^*\mathrm{d}\theta
    \\[2mm]
    = \int_{\overline{E}^\eps}|\Omega^\eps|^{-1}\left(\int_{\Omega^\eps}\frac{\partial
    u(\theta,s)}{\partial
    \theta}\mathrm{d}\mathbf{s}\right)^2\mathrm{d}\theta
    \leq \!\int_{E^\eps}\!\!\left|\nabla u\right|^2\!\dxx.
\end{multline}
For  $V^\eps$, we have by Lemma \ref{LemmaAB} and (\ref{Pdef})
\begin{multline}\label{eqProofPeps1}
\eps^{(N-1)}\int_{\overline{V}^\eps}
    \left| \left(P^\eps u\right)^{'}\right|^2\rho^*\dtet
    =\eps^{(N-1)} \!\sum_{e, \tilde{e}\in N(v)}u_e\overline{u_{\tilde{(e)}}} \!\!
    \int_{\overline{V}^\eps}
    (\psi_{(e)}^{\eps})'(\psi_{\tilde{(e)}}^{\eps})'
    \rho^*\dtet
\\[2mm]
=\eps^{(N-2)}\delta^{-1}{\vec{u}}{\bf \overline{A}}(\vec{u})^*
\leq(\eps\delta)^{(N-2)}\alpha^{\overline{A}}
|{\vec{u}}\llcorner{\vec{1}}|^2 \\[2mm]
\leq
(\eps\delta)^{(N-2)}\frac{\alpha^{\overline{A}}}{\beta^A}\int_{V}|\nabla
\widehat{u}|^2 \dxx=\frac{\alpha^{\overline{A}}}{\beta^A}
\int_{V^\eps} |\nabla u|^2\dxx  .
\end{multline}


3. In the edges $E^\eps$, we use the same argument as in
\cite{Rub1}. By the inequality
\begin{equation}\label{LemmaPEq1}
    (a+b)^2\geq(1-\sqrt{\eps})a^2-b^2/\sqrt{\eps},
\end{equation}
we have that
\begin{multline*}
\eps^{(N-1)}\int_{\overline{E}^\eps}|P^\eps u|^2\rho^*\dtet
=\int_{E^\eps}|P^\eps u|^2 \dss\dtet=\int_{E^\eps}|u+(P^\eps
u-u)|^2\dss\dtet
\\
\geq
\int_{\overline{E}^\eps}\int_{\Omega^\eps}\left[(1-\sqrt{\eps})|u(\theta,\mathbf{s})|^2-
\frac{1}{\sqrt{\eps}}|P^\eps
u-u(\theta,\mathbf{s})|^2\right]\mathrm{d}\mathbf{s}\mathrm{d}\theta
\\
= \eps^{(N-1)}\int_{\overline{E}^\eps}\int_{\Omega}
\left[(1-\sqrt{\eps})|\widetilde{u}(\theta,\mathbf{s}|^2)
-\frac{1}{\sqrt{\eps}}|P
\widetilde{u}-\widetilde{u}(\theta,\mathbf{s})|^2\right]\mathrm{d}\mathbf{s}\mathrm{d}\theta.
\end{multline*}
Notice that for each $\theta$ we have that
$P\widetilde{u}(\theta)-\widetilde{u}(\theta,\mathbf{s})$ has
average zero on $\Omega$. By Poincar\'{e} inequality in
$H^1(\Omega)$, there exists a constant $D>0$ such that
$\int_{\Omega}|P \widetilde{u}-\widetilde{u}|^2\dss \leq D
\int_{\Omega}|\nabla \widetilde{u}|^2\dss$ and hence,
\begin{multline*}
\eps^{(N-1)}\!\int_{\overline{E}^\eps}\!\!|P^\eps
u|^2\rho^*\mathrm{d}\theta \!\geq\! \eps^{(N-1)}\!
\int_{\overline{E}^\eps}\!\int_{\Omega}\!\!\left[(1\!-\!\sqrt{\eps})
|\widetilde{u}(\theta,\mathbf{s})|^2\!-\!\frac{1}{\sqrt{\eps}}D|\nabla
\widetilde{u}(\theta,\mathbf{s})|^2\right]\mathrm{d}\mathbf{s}\mathrm{d}\theta
\\
 = \int_{E^\eps}\left[(1-\sqrt{\eps}) |u|^2-\eps^{3/2}D|\nabla u|^2\right]\dxx .
\end{multline*}
Therefore,
$$\int_{\cup_eE^\eps(e)}[(1-\sqrt{\eps})|u|^2
-\eps^{3/2}D|\nabla u|^2]\dxx
\leq\eps^{(N-1)}\int_{T_1}\rho^*|P^\eps u|^2\dtet.
$$
On the other hand, by Lemma \ref{LemmaConvergenceAtEnds},
$$\int_{\cup_vV^\eps(v)}\!\!\![(1\!-\!\sqrt{\eps})
|u|^2\!-\!\eps^{3/2}D|\nabla u|^2]\!\dxx\! \leq\!
\int_{\cup_vV^\eps(v)}\!\!\!(1\!-\!\sqrt{\eps})|u|^2\!\dxx\!
\leq\! c\eps(1-\sqrt{\eps})\!\int_{T_N^\eps}\!\!|\nabla u|^2
\!\dxx.
$$
Summing the last two inequalities, we obtain the proof of  part 3.

4. Since  $$\eps^{(N-1)}\!\!\int_{\overline{E}^\eps}
W_{T_1,\eps}|P^\eps u|^2\rho^*\dtet
=\int_{E^\eps}\!W_{T_N^\eps}|P^\eps u|^2\dss\dtet ,
$$  it is sufficient to prove for the edges that
$$\int_{E^\eps}W_{T_N^\eps}|P^\eps
u|^2\dss\dtet \leq\int_{E^\eps}(1+2\sqrt{\eps})
W_{T_N^\eps}|u|^2+O(\eps)|\nabla u|^2) \dxx .$$
Using \eqref{LemmaPEq1}, we have that
$$\int_{E^\eps}W_{T_N^\eps}|u|^2\dxx
\geq (1-\sqrt{\eps})\int_{E^\eps}W_{T_N^\eps}|P^\eps
u|^2\dss\dtet-\frac{1}{\sqrt{\eps}}
\int_{E^\eps}W_{T_N^\eps}|u-P^\eps u|^2 \dss\dtet .
$$
Therefore, if $0<\eps<1$ is small enough so that $1\leq
(1-\sqrt{\eps})(1+2\sqrt{\eps})$, then by Poincar\'{e} inequality,
there exists a constant $D$ such that
\begin{multline}\label{64}
\int_{E^\eps}\!\!\!\!W_{T_N^\eps}|P^\eps u|^2\dss\dtet
\leq(1+2\sqrt{\eps})\left\{\int_{E^\eps}\!\!\!\!W_{T_N^\eps}|u|^2\dxx
+\frac{1}{\sqrt{\eps}}
\int_{E^\eps}\!\!\!\!W_{T_N^\eps}|u\!\!-\!\!P^\eps
u|^2\right\}\dss\dtet
\\[2mm]
\leq(1+2\sqrt{\eps})\int_{E^\eps}W_{T_N^\eps}|u|^2\dxx
+C_W(1+2\sqrt{\eps})\frac{1}{\sqrt{\eps}}
\int_{\overline{E}^\eps}\int_{\Omega}|\widetilde{u}-P
\widetilde{u}|^2(\eps\delta)^{(N-1)}\dss\dtet
\\[2mm]
\leq (1+2\sqrt{\eps})\int_{E^\eps}W_{T_N^\eps}|u|^2\dxx
+C_WD(1+2\sqrt{\eps})\frac{(\eps\delta)^2}{\sqrt{\eps}}
\int_{E^\eps}|\nabla_\mathbf{s}u|^2\dss\dtet
\\[2mm]
\leq(1+2\sqrt{\eps})\int_{E^\eps}W_{T_N^\eps}|u|^2\dxx
+O(\eps^{3/2}) \int_{E^\eps}|\nabla u|^2\dxx .
\end{multline}
For the connectors we obtain  by Lemma
\ref{LemmaL2NormNotInEndsT1} and part 2,
$$\eps^{(N-1)}\int_{\cup_v\overline{V}^\eps(v)}
W_{T_1,\eps}|P^\eps u|^2\rho^* \dtet\leq \eps^{N}C_W\int_{T_1}
\left|\left( P^\eps u\right)^{'}\right|^2\rho^* \dtet\leq \eps
C_W\int_{T_N^\eps}|\nabla u|^2 \dxx$$ which, together with
(\ref{64}), yields the proof of part 4.
\end{proof}
\begin{corollary}
\label{P_eps_Theorem} For all $\eps>0$ sufficiently small, there
exists a constant $c>0$ such that the Rayleigh quotients
$$R_{H^1_{0,\rho^*}(T_1)}[P^{\eps}u]
:= \dfrac{\int_{T_1} \left(|(P^{\eps}
u)^{'}|^2\rho_P^\eps+W_{T_1,\eps}|P^\eps
u|^2\rho^*\right)\mathrm{d}\theta
}{\int_{T_1}|u|^2\rho^*\mathrm{d}\theta}\,,
$$
and
$$
R_{H_0^1(T_N^\eps)}[u]
:=\frac{\int_{T_N^\eps}(|\nabla
u|^2+W_{T_{n}^{\eps}}|u|^2)\mathrm{d}\mathbf{x}
}{\int_{T_N^\eps}|u|^2\mathrm{d}\mathbf{x}}
$$
satisfy
\begin{equation}\label{ReqSatP}
    R_{H^1_{0,\rho^*}(T_1)}[P^{\eps}u]\leq
    \frac{[1+O(\sqrt{\eps})]R_{H_0^1(T_N^\eps)}[u]}
    {1-\sqrt{\eps}-c\eps R_{H_0^1(T_N^\eps)}[u]}
   \qquad  \forall u\in H_0^1(T_N^\eps).
\end{equation}
\end{corollary}

\begin{remark} {\em Notice that
$R_{H^1_{0,\rho^*}(T_1)}[P^{\eps}u]$ is the Rayleigh quotient of
the width-weighted operator $A_\eps$ defined in
\eqref{DefWidthWeightedA}, substituting
$\rho_{\alpha,\eps}=\rho_P^\eps$ and $\rho_{b,\eps}=\rho^*$. }
\end{remark}

\subsection{$T_1$-based estimates for the spectrum on
$T^\eps_N$} \label{SectionKobyLemma}

Rubinstein and Schatzman have proved the following general lemma
\cite{Rub1}.
\begin{lemma} \label{Lemma5Koby}
Let $A_j$ be bounded below, selfadjoint operators defined on
Hilbert  spaces $H_j$, where $j=0,1$, and let $\{\lambda_m(A_j)\}$
be the nondecreasing sequence of the corresponding eigenvalues.
    Denote by $D_j$ the domain of the maximal quadratic form
    associated with $A_j$ and by $R_j$ the Rayleigh quotient
    associated with $A_j$. Suppose that there exists a continuous
    linear operator $S$ mapping $D_1$ to $D_0$ and an increasing function
    $\phi:\mathbb{R}\to \mathbb{R}\cup \{+\infty\}$ such
    that $\exp(-\phi)$ is continuous, and
    $$ \>\>\>R_0(Su)\leq
    \phi(R_1(u))\qquad \forall u\in D_1\setminus \ker(S).$$
    Assume that for a given $m$,
    \begin{equation}\label{kobykon}\mu:=\inf\{R_1(v)\mid v\!\in \!D_1\cap \ker (S),v\neq
    0\}>\lambda_m(A_1).\end{equation}
    Then
    \begin{equation} \label{RubEq}
        \lambda_m(A_0)\leq \phi(\lambda_m(A_1)).
    \end{equation}
\end{lemma}
Using Lemma \ref{Lemma5Koby}, we obtain bounds for the eigenvalues
of $T_N^\eps$.  Let $\nu_m^\eps$ denotes the $m$-th eigenvalue of
the Schr\"{o}dinger operator
$$L_\eps:=-\Delta+W_{T_N^\eps} \ . $$ Denote the operators
$$A^\eps_Q:= -\frac{1}{\rho^*}
    \frac{\mathrm{d}}{\mathrm{d}\theta}\left(\rho_Q^\eps
    \frac{\mathrm{d}}{\mathrm{d}\theta}\right)+W_{T_1,\eps}\,,
     \qquad A^\eps_P:=-\frac{1}{\rho^*}
    \frac{\mathrm{d}}{\mathrm{d}\theta}\left(\rho_P^\eps
    \frac{\mathrm{d}}{\mathrm{d}\theta}\right)+W_{T_1,\eps}  \ , $$
and let $\mu_m^\eps$ (resp. $\lambda_m^\eps$) be the $m$-th
eigenvalue of $A_Q^\eps$ (resp. $A_P^\eps$). We will omit the
superscript $\eps$  in $\nu_m^\eps$, $\mu_m^\eps$, and
$\lambda_m^\eps$ whenever there is no danger of confusion.
\begin{theorem}\label{TheoremEigenvaluesBounds}
    Using the notations above, for all $M\in \mathbb{N}$ there exist $\eps_M>0$
and a constant $c>0$ such that for all $m\leq M$  and
$0<\eps<\eps_M$,  we have
    \begin{equation}\label{TheoremEigenvaluesBoundsEq1}
         \nu_m^\eps\leq \phi_Q^\eps(\mu_m^\eps),
    \end{equation}
    and
    \begin{equation}\label{TheoremEigenvaluesBoundsEq2}
        \lambda_m^\eps \leq \phi_P^\eps(\nu_m^\eps),
    \end{equation}
     where
    $$\phi_Q^\eps(x):=\!\left\{\!\!%
\begin{array}{ll}
    \dfrac{(1\!+\!c\eps)x}{1\!-\!c\eps x} & x\!<\! {(c\eps)^{-1}}, \\[3mm]
    +\infty & \text{otherwise},
\end{array}%
\right.
    \>\>\>\>\mbox{ and }\>\>\>\>
    \phi_P^\eps(x):=\!\left\{\!\!
    \begin{array}{cc}
      \dfrac{(1\!+\!c\eps)x}{1\!-\!\sqrt{\eps}\!-\!
      c\eps x}
      & x\!<\! \dfrac{1\!-\!\sqrt{\eps}}{c\eps}\,, \\[3mm]
      +\infty & \text{otherwise}.
    \end{array}
    \right.$$
  \end{theorem}
\begin{proof}[{Proof of Theorem~\ref{TheoremEigenvaluesBounds}}]

Without loss of generality, we assume  that $W_{T_1}$ is positive.
In order to prove \eqref{TheoremEigenvaluesBoundsEq1}, we wish to
apply Lemma \ref{Lemma5Koby} on $S\!=\!Q^{\eps},
D_0\!=\!H^1_0(T_N^\eps)$, $ D_1\!=\!H^1_{0,\rho^*}(T_1),
A_0\!=\!L_\eps, A_1\!=\!A^\eps_Q$ , $R_0\!=\!R_{H^1_0(T_N^\eps)}$,
and $ R_1\!=\!R_{H^1_{0,\rho^*}(T_1)}$.
We, therefore,  show that there exists $C>0$ such that for any
$\eps>0$
\begin{equation}\label{ReqKerQ}
   \inf\left\{R_{H^1_{0,\rho^*}(T_1)}[f]\mid f\in \ker Q^\eps, \> f\neq 0\right\} \geq
   \frac{1}{C\eps^2}\,.
\end{equation}
  Indeed,
$$\ker Q^{\eps}
 =\left\{f\in H_{0,\rho^*}^1(T_1)\mid f(\theta)=0
 \>\>\>\>\>\forall\theta\!\in\! T_1\backslash
 \cup_{v}\overline{V}^\eps \right\}.$$ Therefore, in order to estimate
$R_{H^1_{0,\rho^*}(T_1)}[f]$ for $f\in \ker(Q^\eps)$, we actually
need to estimate this quotient in each component
$\overline{V}^\eps\cap e$. However, we have that
 $$|f(\theta)|^2
  =\left|\int_{p^\eps}^{\theta}f^{'}\mathrm{d}\vartheta\right|^2
 \leq
 |p^\eps-\theta|\int_{p^\eps}^{\theta}|f^{'}|^2\mathrm{d}\vartheta,
$$
where $p^\eps\in\partial\overline{V}^\eps$. Multiply the above by
$\rho^*$ (which, we recall, is constant on each component
$\overline{V}^\eps\cap e$), we find the existence of
 $C>0$ such that
$$\int_{\overline{V}^\eps\cap e}|f|^2\rho^*\mathrm{d}\theta
\leq \int_{\overline{V}^\eps\cap e}\left(
|p^\eps-\theta|\int_{p^\eps}^{\theta}|f^{'}|^2\rho^*
\mathrm{d}\vartheta \right)\mathrm{d}\theta \leq C\eps^2
\int_{\overline{V}\cap e}|f^{'}|^2\rho^*\dtet .
$$
Thus, (\ref{ReqKerQ}) is verified provided $\eps$ is sufficiently
large. Hence, \eqref{TheoremEigenvaluesBoundsEq1} follows from
Corollary \ref{Q_eps_Theorem} and Lemma \ref{Lemma5Koby}.

In order to prove \eqref{TheoremEigenvaluesBoundsEq2}, we wish to
apply Lemma \ref{Lemma5Koby} to $S\!=\!P^{\eps},
D_0\!=\!H^1_{0,\rho^*}(T_1)$, $ D_1=H^1_0(T_N^\eps), A_0=A_P^\eps,
A_1=L_\eps, R_0=R_{H^1_{0,\rho^*}(T_1)},\>\>\mbox{ and }\>\>
R_1=R_{H^1_0(T_N^\eps)}$.
To this end,  we show that there exists $C>0$ such that for any
$\eps>0$
\begin{equation}\label{ReqKerP}
   \inf\left\{R_{H_0^1(T_N^\eps)}[u]\>\>|\>\> u\in \ker P^\eps, \> u\neq 0\right\} \geq
   \frac{1}{C\eps}.
\end{equation}
We notice that if $u\in \ker P^\eps$, then its averages on the
cross sections $\Omega_j^\eps$ of $E_j^\eps$  vanish. Therefore,
using the $(N-1)$-dimensional  Poincar\'{e} inequality for
functions whose average is  zero, we obtain that there is a
constant $D$ such that:
\begin{equation}\label{771}\int_{E^\eps}\!\!|u|^2\dss\dtet
\leq D\eps^{2(N-1)}\!\!\int_{\overline{E}^\eps}\!\!
\int_{\Omega^\eps}\!\!|\nabla_su|^2\dss\dtet \leq
D\eps^{2(N-1)}\int_{E^\eps}\!\!\! (|\nabla
u|^2+W_{T_N^\eps}|u|^2)\dss\dtet .
\end{equation}
By Lemma \ref{LemmaConvergenceAtEnds}, there exists $C>0$ such
that for any  $u\in H_0^1(T_N^\eps)$
\begin{equation}\label{772}
\int_{\cup_{v}V^\eps(v)}|u|^2\dxx \leq
C\eps\int_{T_N^\eps}(|\nabla u|^2+W_{T_N^\eps}|u|^2) \dxx .
\end{equation} Therefore, (\ref{771}) and (\ref{772}) imply
\eqref{ReqKerP}. Thus, \eqref{TheoremEigenvaluesBoundsEq2} follows
by Corollary \ref{P_eps_Theorem} and Lemma \ref{Lemma5Koby}.
\end{proof}
\begin{remark} {\em Theorem~\ref{TheoremEigenvaluesBounds} is similar to
\cite[Theorem5]{Rub1} proved for a finite graph with a
constant-width thin domain. }
\end{remark}
\begin{theorem}
  For each $m\in\mathbb{N}$, the $m$-the eigenvalue of  the Schr\"{o}dinger operator $L_\eps$ on
    $H^1_0(T_N^\eps)$ converges as $\eps\rightarrow 0$
    to the $m$-the eigenvalue of limit width-weighted operator $A$ on $H_{0}^1(T_1)$.\end{theorem}
\begin{proof}
    We use in this proof the notations of
    Theorem~\ref{TheoremEigenvaluesBounds}. Notice that for small enough $\eps$,
    $\phi_Q^\eps$ and $\phi_P^\eps$ are continuous monotone
    increasing function, which satisfy $$\lim_{\eps\rightarrow
    0}\phi_Q^\eps(x) =x,\>\>\>
    \lim_{\eps\rightarrow 0} \phi_P^\eps (x)=x.$$
    Moreover, since the operators we refer to in
    Theorem~\ref{TheoremEigenvaluesBounds} satisfy the conditions of
    Theorem~\ref{MainTheorem}, we have for each $m\in \mathbb{N}$ that both $\mu_m^\eps$
    and $\lambda_m^\eps$ (see \eqref{TheoremEigenvaluesBoundsEq1} and
    \eqref{TheoremEigenvaluesBoundsEq2}) converge as $\eps\rightarrow 0$ to
    the $m$-th eigenvalue of the limit width-weighted operator $A$.
    Since $A$ has a discrete spectrum, the result follows.
\end{proof}
\mysection{Convergence of eigenfunctions of Laplace operator on
$T_N^\eps$} \label{SecEigHarmonicConv} In \cite{Kosugi1,Kosugi2},
Kosugi has proved that the solutions of $\Delta u+f(u)=0$ in thin
network-shaped bounded domains that satisfy Neumann boundary
condition, converge to solutions of appropriate equations on the
skeleton of the domain. In \cite{Kosugi1}, Kosugi deals only with
domains which are formed by joining straight tubes around some
graph, while in \cite{Kosugi2} the results are extended to general
domains around graphs. However, trees with infinite number of
vertices and nonsmooth boundaries are not considered in these
papers. Using the transformation $P^\eps$ developed for Theorem
\ref{P_eps_Theorem}, we give a simple proof for the convergence of
projections into $H^1_{0,*}(T_1)$ of eigenfunctions $u_\eps$ of
the  Laplace operator on $H_0^1(T_N^\eps)$. Specifically, we show
in Theorem~\ref{TheoremConvergenceEigenfun} that $P^\eps u_\eps$
converges to eigenfunctions of the following limit width-weighted
operator on $T_1$
$$ L_* u:= \left(\rho^*\right)^{-1}\left( \rho^* u^{'}\right)^{'}
. $$

 First, we need to prove the following auxiliary Lemma.
\begin{lemma} \label{LemmaConvergenceToCont}
    Assume that $u\in H_0^1(T_n^\varepsilon)$ satisfies
    $||u||_{H^1_0(T_n^\varepsilon)}=\varepsilon^{(n-1)/2}.$
    Fix a vertex $v$, and denote by $p_e$ the `end point' in
    $\overline{V}^\varepsilon\cap \overline{E}^\varepsilon(e)$.
    Then there is a constant $C$ which depends on $v$ but is independent on $\varepsilon$
    such that for $e,\tilde{e}\in N(v)$ we have
        \begin{equation} \label{ContPEq}
            |P^\varepsilon u(p_e)-P^\varepsilon u(p_{\tilde{e}})|
            \leq C\sqrt{\mathrm{dist}(p_e,p_{\tilde{e}})}\,,
        \end{equation}
    where $\mathrm{dist}(\cdot,\cdot)$ is the standard distance function on $T_1.$
    \end{lemma}
{\bf Proof.}  Notice that since $C^1(T_n^\varepsilon)$ is dense in
$H_0^1(T_n^\varepsilon)$ we may assume without loss of generality
that $u\in C^1(T_n^\varepsilon).$\\
Let $q,r\in \overline{V}^\varepsilon\cap e.$ By
\eqref{eqProofPeps1},
\begin{multline*}
|P^\varepsilon u(q)-P^\varepsilon u(r)|^2
\!=\!\left|\int_{r}^{q}\frac{\mathrm{d}}{\mathrm{d}\theta}(P^\varepsilon
u)\mathrm{d}\theta\right|^2
 \!\!\leq
\mathrm{dist}(q,r)\frac{1}{\rho^{*}_e}\int_{r}^{q}
\!\left|\frac{\mathrm{d}}{\mathrm{d}\theta}(P^\varepsilon
u)\right|^2\!\!\rho^*\mathrm{d}\theta\\[2mm]
 \leq
\mathrm{dist}(q,r)\frac{\varepsilon^{1-n}}{\rho^{*}_e}\frac{\alpha^{\overline{A}}}{\beta^A}
\int_{V^\varepsilon}|\nabla
u|^2\mathrm{d}\mathbf{s}\mathrm{d}\theta \leq
\mathrm{dist}(q,r)\frac{\varepsilon^{1-n}}{\rho^{*}_e}\frac{\alpha^{\overline{A}}}{\beta^A}
\varepsilon^{n-1} \leq C\mathrm{dist}(q,r)
\end{multline*}
for some constant $C.$ Therefore, $|P^\varepsilon
u(p_e)-P^\varepsilon u(p_{\tilde{e}})|\leq
 2C\sqrt{\mathrm{dist}(p_e,p_{\tilde{e}})}$.
\qed

\begin{theorem} \label{TheoremConvergenceEigenfun}
    Let $u_\eps\in H_0^1(T_N^\eps)$ be an eigenfunction with eigenvalue
    $\lambda_\eps$ of the Laplace operator on $T_N^\eps$,
    such that $\|u_\eps\|_{L^2(T_N^\eps)}=\eps^{(N-1)/2}$.
    Assume that $\lim_{\eps\rightarrow
    0}\lambda_\eps=\lambda^*. $
    Then there exists an eigenfunction $u^*$ of $L_*$  which corresponds to $\lambda^*$,
    such that up to a subsequence,
    $$ u^*=\lim_{\eps\rightarrow 0}P^\eps
    u_\eps$$ locally uniformly.
\end{theorem}
\begin{proof} By elliptic regularity, $u_\varepsilon\in C^2(T_N^\eps).$
Our proof consists of three steps.

\noindent {\bf Step 1.} Let us show that $P^\eps u_\eps$ converges
to a solution $u^*$ of $\frac{d^2u}{d\theta^2}=\lambda^* u$ on
each edge of $e\in T_1$.

By parts 2 and 4  of Lemma~\ref{LemmaP} (with $W= 1$), we obtain
that $P^\eps u_\eps$ are uniformly bounded in $ H_{*,1}(T_1)$.
This implies, in particular, that $P^\eps u_\eps$ are uniformly
locally bounded in $L^\infty(T_1)$. In addition, (up to a
subsequence) $\lim_{\eps\rightarrow 0} P^\eps u_\eps=u^*$ holds
locally uniformly by Arzel\`{a}-Ascoli's Theorem. Fix an edge
$e\in T_1$, and $\theta_1,\theta_2\in e$. Let $\zeta(\theta)\in
C_0^\infty([\theta_1,\theta_2])$.  If $\eps>0$ is sufficiently
small, then $\theta_1,\theta_2\in \overline{E}^\eps$. Therefore,
\begin{multline*}
\int_{\theta_1}^{\theta_2}P^\varepsilon u_\varepsilon(\theta)
\zeta''(\theta)d\theta
=\int_{\theta_1}^{\theta_2}\frac{1}{|\Omega^\varepsilon|}
\left(\int_{\Omega^\varepsilon} u_\varepsilon
(\theta,\mathbf{s})d\mathbf{s}\right) \zeta''(\theta)\dtet
 \\
=\frac{1}{|\Omega^\varepsilon|}
\int_{\Omega^\varepsilon}\int_{\theta_1}^{\theta_2} u_\varepsilon
(\theta,\mathbf{s}) \Delta\zeta(\theta)\dtet \dss
=-\frac{1}{|\Omega^\varepsilon|}
\int_{\Omega^\varepsilon}\int_{\theta_1}^{\theta_2} \nabla
u_\varepsilon (\theta,\mathbf{s})\cdot \nabla\zeta(\theta)\dtet \dss \\
=-\frac{\lambda_\varepsilon}{|\Omega^\varepsilon|}
\int_{\Omega^\varepsilon}\int_{\theta_1}^{\theta_2} u_\varepsilon
(\theta,\mathbf{s})\zeta(\theta)\dtet\dss
=-\lambda_\varepsilon \int_{\theta_1}^{\theta_2} P^\varepsilon
u_\varepsilon (\theta)\zeta(\theta)\dtet.
\end{multline*}
Hence, $P^\varepsilon u_\varepsilon\in H^2([\theta_1,\theta_2])$
and $-(P^\varepsilon u_\varepsilon)''=\lambda_\varepsilon
P^\varepsilon u_\varepsilon$ in the weak sense and by elliptic
regularity also in the strong sense. Moreover $P^\varepsilon
u_\varepsilon$ is $C^\infty$ in $\overline{E}^\eps$. Since
$\lambda_\eps\rightarrow\lambda^*$ and $P^\eps u_\eps\rightarrow
u^*$ uniformly on $e$, the second derivatives $(P^\varepsilon
u_\varepsilon)''$ converge uniformly to $(u^*)''$, which also
implies the same convergence for the first derivatives
$(P^\varepsilon u_\varepsilon)'$.
\par\noindent  {\bf Step 2.} We show now that $u^*$ is in the domain
of $L_*$.  For this, we must only show that $u^*$ satisfies the
corresponding  Kirchhoff's conditions. The continuity at the
vertices is satisfied by Lemma \ref{LemmaConvergenceToCont}. The
second Kirchhoff condition is given by
$$\sum_{e\in N(v)} \rho_e^* u_e'(v)=0,$$
where $N(v)$ is the set of all edges adjacent to the vertex $v$.
Recall that $\rho_e^*=\delta_e^{(N-1)}|\Omega_e|$ takes a constant
value on each edge $e$.

 Let
$\overline{U}\subset T_1$ be a neighborhood of the vertex $v$
which contains no other vertex, and let $\theta_e\in
\partial\overline{U}$ be the point of $\partial\overline{U}$
contained in $e\in N(v)$. Let $U^\eps\subset T_N^\eps$ be the {\it
inflation} of $\overline{U}$, that is, $\overline{U}=U^\eps \cap
T_1$.
 In particular, for sufficiently small $\eps$ we have
$$\partial U^\eps=\left(U^\eps\cap\partial T^\eps_N\right)
\bigcup_{e\in N(v)}S_e\,,$$ where $S_e= \{ \mathbf{s}; (\theta_e,
\mathbf{s})\in E^\eps(e)\}$.
Let $\zeta_\varepsilon\in C^\infty(U^\varepsilon)$ be a function
which does not depend on $\mathbf{s}$ in the edges, satisfies
$\zeta_\varepsilon(\mathbf{x})=1$ for all $x\in V^\varepsilon(v)$,
$0\leq \zeta_\varepsilon(\mathbf{x})\leq 1$ for all $x\in
U^\varepsilon$, and vanishes around each $S_e$. Since
$u_{\varepsilon}$ is an eigenfunction, we have
\begin{multline*}
\lambda_{\varepsilon} \int_{U^\varepsilon}u_\varepsilon
\zeta_\varepsilon \dxx
=\int_{V^\varepsilon(v)}\nabla  u_\varepsilon \cdot\nabla
\zeta_\varepsilon \dxx +\int_{U^\varepsilon(v)\backslash
V^\varepsilon(v)}\nabla u_\varepsilon\cdot \nabla
\zeta_\varepsilon \dxx
\\
=\int_{U^\varepsilon(v)\backslash V^\varepsilon(v)} \frac{\partial
u_\varepsilon}{\partial \theta} \frac{
\mathrm{d}\zeta_\varepsilon}{\mathrm{d}\theta} \dss\dtet.
\end{multline*}
As $\zeta_\varepsilon$ depends only on $\theta$ on
$U^\varepsilon(v)\backslash V^\varepsilon(v)$ and equals one at
$p_e$, we get
\begin{multline}\label{xzx1}
\int_{U^\varepsilon(v)\backslash V^\varepsilon(v)} \frac{\partial
u_\varepsilon}{\partial \theta} \frac{
\mathrm{d}\zeta_\varepsilon}{\mathrm{d}\theta}
\mathrm{d}\mathbf{s}\mathrm{d}\theta = \sum_{e\in
N(v)}|\Omega_e^\varepsilon|
 \int_{p_e}^{\theta_e} \frac{\partial P^{\varepsilon}u_\varepsilon}{\partial \theta}
\frac{ \mathrm{d}\zeta_\varepsilon}{\mathrm{d}\theta}\dtet \\
 = - \sum_{e\in N(v)}|\Omega_e^\varepsilon|\left[\frac{\partial
P^{\varepsilon}u_\varepsilon}{\partial
\theta}(p_e)\zeta_\varepsilon(p_e) + \int_{p_e}^{\theta_e}
(P^{\varepsilon} u_{\varepsilon})'' \zeta_\varepsilon
\dtet\right] \\
=  -\sum_{e\in N(v)}|\Omega_e^\varepsilon|\left[\frac{\partial
P^{\varepsilon}u_\varepsilon}{\partial \theta}(p_e) -
\lambda_\varepsilon \int_{p_e}^{\theta_e}
P^{\varepsilon}u_{\varepsilon} \zeta_\varepsilon \dtet\right]
 \\
 = - \sum_{e\in N(v)}|\Omega_e^\varepsilon|\frac{\partial
P^{\varepsilon}u_\varepsilon}{\partial \theta}(p_e) +
\lambda_\varepsilon \int_{U^\varepsilon(v)\backslash
V^\varepsilon(v)} u_{\varepsilon} \zeta_\varepsilon \dxx.
\end{multline}
The change of order of integration and differentiation in the
first line of \eqref{xzx1} is easily justified by approximating
$u_\varepsilon$ with a smooth function. We therefore obtain that
$$
\sum_{e\in N(v)}|\Omega_e^\varepsilon|\frac{\partial
P^{\varepsilon}u_\varepsilon}{\partial \theta}(p_e)
=-\lambda_{\varepsilon} \int_{V^\varepsilon(v)}u_\varepsilon
\zeta_\varepsilon \dxx,
$$
and since $|V^\varepsilon(v)|=c\varepsilon^{N},$ we arrive at the
estimate
\begin{align}\label{xzx2}
\left|\sum_{e\in N(v)}\rho^{*}_e\frac{\partial
P^{\varepsilon}u_\varepsilon}{\partial \theta}(p_e)\right|
\leq &\
c\varepsilon^{1-N}\varepsilon^{N/2}\lambda_{\varepsilon}\left(
\int_{V^\varepsilon(v)}u_\varepsilon^2(\mathbf{x})
\dxx\right)^{1/2}
\\[0.2cm]
= &\ c\lambda_\varepsilon \varepsilon^{1/2}.\notag
\end{align}
Letting $\eps\to 0$, we obtain by Step 1 that the left hand side
of \eqref{xzx2} converges to
$$\left|\sum_{e\in
N(v)}\rho^*_{e}(u_e^*)'(v)\right|$$ and the right hand side to
zero.
\par\noindent {\bf Step 3}. It remains to prove that $u^*\not\equiv
0$. Let $T_{N,j }$ denote the $j$ first generations in $T_N$. By
lemmas \ref{LemmaP} and \ref{LemmaL2NormNotInEndsT1} there are
constants $c,C>0$ and a function $R(j)$ which tends to zero as
$j\rightarrow \infty$ such that
\begin{multline*}
 \eps^{(N-1)}\int_{T_{1,j}}
|P^\eps u_\eps|^2\rho^*\mathrm{d}\theta\\
 =\eps^{(N-1)}\int_{T_1}
|P^\eps u_\eps|^2\rho^*\mathrm{d}\theta-
\eps^{(N-1)}\int_{T_1\backslash T_{1,j}} |P^\eps
u_\eps|^2\rho^*\mathrm{d}\theta \\
\geq(1-\sqrt{\eps})\eps^{(N-1)}-c\lambda_\eps\eps^{N}
-\frac{c^2}{C^2}R(j)^2\lambda_\eps\eps^{(N-1)}
\\
=\eps^{(N-1)}\left[(1-\sqrt{\eps})-c\lambda_\eps\eps
-\frac{c^2}{C^2}R(j)^2\lambda_\eps\right].
\end{multline*}
Choose $\eps>0$ small enough and $j$ large enough so that $
(1\!-\!\sqrt{\eps})-c\lambda_\eps\eps
-\frac{c^2}{C^2}R(j)^2\lambda_\eps \!\!\geq\!\!\gamma $ for a
constant $\gamma>0$. Then $ \int_{T_{1,j}} |P^\eps
u_\eps|^2\rho^*\mathrm{d}\theta \geq \gamma>0$. By the local
uniform convergence of $P^\eps u_\eps$ to $u^*$ we have that
$$ \int_{T_1} |u^*|^2\rho^*\dtet
\geq \int_{T_{1,j}} |u^*|^2\rho^*\dtet =\lim_{\eps \rightarrow 0 }
\int_{T_{1,j}} |P^\eps u_\eps|^2\rho^*\dtet \geq \gamma
>0,$$
so, $u^*\not\equiv 0$.
\end{proof}
\begin{center}
{\bf Acknowledgments}
\end{center}
 The paper is based on part of the Ph.~D. thesis
\cite{Z} of Daphne Zelig, completed in 2005 at the Technion, under
the supervision of Moshe Israeli, Yehuda Pinchover and Gershon
Wolansky. The authors would like to thank Professors Peter
Kuchment, Alexander Sobolev, and Michael Solomyak for valuable
discussions.

This work was partially supported by the RTN network ``Nonlinear
Partial Differential Equations Describing Front Propagation and
Other Singular Phenomena", HPRN-CT-2002-00274. The works of Y.~P.
and G.~W. were also partially supported by the Israel Science
Foundation (grants 1136/04 and 406/05, respect.) founded by the
Israeli Academy of Sciences and Humanities, and by the Fund for
the Promotion of Research at the Technion.


%
\end{document}